\documentclass[graybox,english,11pt]{article}
\usepackage{geometry}
\usepackage[latin9]{inputenc}
\usepackage{array}
\usepackage{verbatim}
\usepackage{float}
\usepackage{bm}
\usepackage{multirow}
\usepackage{amsthm}
\usepackage{amsmath}
\usepackage{amssymb}
\usepackage{graphicx}
\usepackage{color}
\usepackage{appendix}
\usepackage{subfig}

\makeatletter

\theoremstyle{plain}
\newtheorem{thm}{\protect\theoremname}[section]
 \theoremstyle{definition}
  
  \theoremstyle{plain}
  \newtheorem{assumption}[thm]{\protect\assumptionname}
  \theoremstyle{plain}
  
  \newtheorem{definition}[thm]{\protect\definitionname}
  \theoremstyle{remark}
  \newtheorem*{rem*}{\protect\remarkname}
  \theoremstyle{plain}
  
  \theoremstyle{remark}
  \newtheorem{rem}[thm]{\protect\remarkname}
  \theoremstyle{plain}
  \newtheorem{algorithm}[thm]{\protect\algorithmname}
 \theoremstyle{plain}
  \newtheorem*{theorem*}{\protect\theoremname}
\theoremstyle{plain}
  \newtheorem{notation}[thm]{\protect\notationname}
\makeatother

\usepackage[english]{babel}
  \providecommand{\algorithmname}{Algorithm}
  \providecommand{\assumptionname}{Assumption}
  \providecommand{\lemmaname}{Lemma}
  \providecommand{\problemname}{Problem}
  \providecommand{\propositionname}{Proposition}
  \providecommand{\remarkname}{Remark}
\providecommand{\theoremname}{Theorem}
\providecommand{\definitionname}{Definition}
\providecommand{\notationname}{Notation}

\usepackage[bottom]{footmisc}
\usepackage{graphics}

\newcommand{\ep}{\epsilon}
\newcommand{\ei}{\epsilon^{-1}}

\newcommand{\R}{\mathbb{R}}

\begin{document}
\title{Parareal multiscale methods for highly oscillatory dynamical systems}
\author{Gil Ariel\thanks{Department of Mathematics, Bar-Ilan University,
Ramat-Gan 5290002, Israel. (arielg@math.biu.ac.il).}
\and Seong Jun Kim\thanks{Department of Mathematics, Georgia Institute of Technology,
Atlanta, GA 30332, USA. (skim396@math.gatech.edu).}
\and Richard Tsai\thanks{Department of Mathematics and Institute for Computational
Engineering and Sciences (ICES), The University of Texas at Austin, TX 78712, USA,  and
KTH Royal Institute of Technology, Sweden. (ytsai@ices.utexas.edu)}}

\date{\date{}}
\maketitle
\begin{abstract}
    We introduce a new strategy for coupling the parallel in time (parareal) iterative methodology with multiscale integrators. 
    Following the parareal framework, the algorithm computes a low-cost approximation of all slow variables
    in the system using an appropriate multiscale integrator, which is refined using parallel fine scale integrations. 
    Convergence is obtained using an alignment algorithm for fast phase-like variables.
    The method may be used either to enhance the accuracy and range of applicability
    of the multiscale method in approximating only the slow variables, or to
    resolve all the state variables.
    The numerical scheme does not require that the system is split into slow and fast coordinates.
    Moreover, the dynamics may involve hidden slow variables, for example, due to resonances. 
    We propose an alignment algorithm for almost-periodic solution, in which case
    convergence of the parareal iterations is proved.
    The applicability of the method is demonstrated in numerical examples.
\end{abstract}

\section{Introduction}
\label{sec:intro} 

The parallel in time, also known as the \lq\lq{}parareal\rq\rq{} method, introduced by Lions, Maday and Turinici \cite{parareal-LMT01} is a simple yet effective scheme for the parallelization of numerical solutions 
for a large class of time dependent problems \cite{parareal-Maday10}.
It consists of a fixed point iteration involving a coarse-but-cheap and a fine-but-expensive integrators.
Computational time is reduced by parallelization of the fine integrations. 
For problems with separated multiple scales,  it is tempting to apply a multiscale solver as a coarse integrator. 
So far, such types of parallel methods are limited to a few special multiscale cases such as chemical kinetics \cite{parareal-Blouze10,parareal-Engblom09,parareal-He10},
dissipative ordinary differential equations (ODEs) \cite{parareal-Legoll13} 
and highly oscillatory (HiOsc) problems in which the oscillatory behavior is relatively simple \cite{parareal-MadayLegoll10,parareal-HS-PDEs13}.
One difficulty stems out from a fundamental difference between the parareal and the multiscale philosophies --- while the former requires point-wise convergence of the numerical solvers (in the state variable), most multiscale schemes gain efficiency by only approximating a reduced set of slowly varying coarse/slow/macroscopic variables \cite{Artstein-Kevrekidis-Slemrod-Titi07,Artstein-Linshiz-Titi07,HMM-review07,GivonKupfermanStuart-review,Equation-free-review09,flavors}.

In this paper, we develop  a general strategy
that couples multiscale integrators and fully resolved fine scale integration for parallel in time computation of HiOsc solutions of a class of ODEs.
There are several advantages in such coupling strategies.
First, some multiscale methods (such as the Poincar\'e-map technique \cite{BFHMM2012})
only approximate the slow constituents or slow variables of the dynamics.
Proper coupling of multiscale and fine scale solvers via a parareal-like framework 
can be efficient (by parallelization) in computing full detailed solutions, 
including the fast phase in the HiOsc dynamics. 
The choice of multiscale method is not limited to the Poincar\'e-map technique. Any multiscale method can be used as a coarse integrator as long as fast phase-like variables are appropriately aligned as required. Then, convergence of the parareal iterations can be shown in a similar manner.
Second, the parareal iterations enhance the stability and accuracy of 
the multiscale scheme, 
in particular when the scale separation in the system is not significant  and
 the corresponding  sampling/averaging errors  are non-negligible.
Finally,  parareal multiscale coupling schemes can deal with more challenging situations,
for example,
(a) \emph{the effective equation is valid almost everywhere macroscopically, 
but is not an adequate description of the system at small but a priori "unpredictable" locations in the phase space
(as these regions may depend on the solutions)}; and (b) \emph{the influence of microscopic solutions in these regions on the 
macroscopic solution elsewhere is significant.} 

In \cite{parareal-Legoll13}, Legoll et at suggest a multiscale parareal scheme
for singularly perturbed ODEs in which the fast dynamics is dissipative, i.e.,
the dynamics relaxes rapidly to a low dimensional manifold.
One of the main contributions of \cite{parareal-Legoll13} is the understanding 
that the slow and fast parts of the dynamics need to be addressed separately.
They suggest two approaches: The first is a straight-forward application of
parareal, which is shown to converge but loses accuracy as the system becomes more singular. 
In Section~\ref{sec:notWorking} we demonstrate that naive parareal does not converge
when applied for HiOsc systems.
The second approach assumes that the system is split into slow and fast variables, or that
a change of variables splitting the system is given.
This approach may be applied to HiOsc systems, but it is relatively restrictive as
in many examples and applications such a splitting is not known. 
Dai et al \cite{parareal-MadayLegoll10} suggest an application of the parareal framework to Hamiltonian systems. They consider two main approaches. The first is a time-reversible iteration scheme (applied together with time-reversible fine and coarse integrators). The second projects solutions at coarse time segments onto the constant energy manifold. The two approaches are also combined together. The first approach is specific to Hamiltonian dynamics and not to general HiOsc problems. The projection method cannot be applied to the HiOsc case because the main difficulty is not with the approximation of slow variables (or constants of motion), but with the fast phase.
In addition, since the methods presented in \cite{parareal-MadayLegoll10} are not multiscale, their accuracy and efficiency are expected to 
deteriorate when the frequencies of oscillations are large.
 Combining the symmetric approach of Dai et al with our alignment method for Hamiltonian systems may be an interesting application, but is beyond the scope of the current manuscript. In particular, the alignment algorithms should also be made symmetric, similar to the ideas of Dai et al.
Applications of parareal methods to Hamiltonian dynamics is also analyzed in  \cite{GanderHairer14}.
Additional approaches to use symplectic integrators with applications to molecular
dynamics include \cite{Audouze-HAL2009,Bal-LNCSE08,Jimenez2011}.
Finally, Haut and Wingate \cite{parareal-HS-PDEs13} suggest a parareal
method for PDEs with linear HiOsc forcing,
As in \cite{parareal-Legoll13}, their method applies exact knowledge of the 
fast variable (the phase in the HiOsc case) to design a convergent parareal scheme.
In this respect, the method proposed here goes further and is also applicable to nonlinear HiOsc systems. However, in this paper the discussion is restricted to the ODE case.
One of the main goals of the current paper is to design a convergent parareal algorithm
that does not require explicit knowledge of the fast and slow variables.

We begin with a short overview of the parareal method within the context of ODEs  and 
test its performance on a simple example HiOsc system.

\subsection{The parareal method for ODEs}
  
Consider the following initial value problem
\begin{equation}
   \dot{u} = f(t,u),~~~u(0)=u_0 ,
\label{eq:generalODE}
\end{equation}
where $u \in \R^d$ and $t\in[0,T]$. 
We assume that $f$ is sufficiently smooth.
Let $H$ denote an intermediate time step, $0<H <T$ and $N=T/H$ an integer.
Suppose that we are given two approximate integrators for \eqref{eq:generalODE}: a cheap coarse integrator with low accuracy denoted ${\mathcal C}$, 
and a fine, high accuracy integrator which is relatively expensive in terms of efficiency, denoted ${\mathcal F}$.
The approximate propagation operators to time $H$ obtained using the the coarse and fine integrators  
are denoted by ${\mathcal C}_H$ and ${\mathcal F}_H$, respectively.

\begin{figure}[tbh]
   \begin{centering}
   \includegraphics[width=11cm]{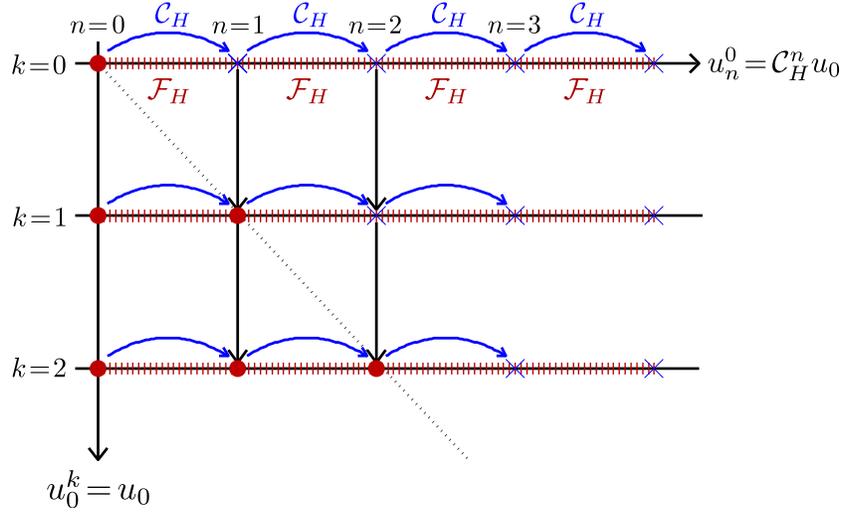}
   \par\end{centering}
   \caption{
   \label{fig:pararealSketch} 
         A sketch depicting the parareal methodology.
         Each parareal iteration is constructed using three integrations: 
         Fine integration starting at $u_n^{k-1}$, 
         coarse integration starting at $u_n^{k-1}$ and another coarse
         integration starting at $u_n^k$. The first two depend only on previous
         iterations and can therefore be computed in parallel.}
\end{figure}

Furthermore, denote by $u^k_n$ the approximation for $u(n H)$ at the $k$'th iteration.
For all iterations, the initial values are the same $\forall k,u^k_0=u_0$.
The objective is to have $u^k_n \to {\mathcal F}_{n H} u_0$ as $k \to \infty$,
i.e., convergence to the approximation given by the high-accuracy fine integrator.
The parareal approximation to \eqref{eq:generalODE}
is as follows.  See Figure~\ref{fig:pararealSketch}  for a sketch of the parareal methodology.
\begin{algorithm}
\label{alg:parareal}
\end{algorithm}
\begin{enumerate}
\item Initialization: Construct the zero'th iteration approximation using a chosen coarse integrator:
   \begin{equation}
       u_0^0 = u_0 ~~ {\rm and} ~~ u^0_{n} = {\mathcal C}_H u^0_{n-1},  ~~~ n=1,\dots,N.
\nonumber
   \end{equation}
\item Iterations: $k=1 \dots K$
   \begin{equation}
       u_0^k = u_0 ~~ {\rm and} ~~ 
      u_{n}^{k} = {\mathcal C}_H u_{n-1}^{k} + {\mathcal F}_H u_{n-1}^{k-1} - {\mathcal C}_H u_{n-1}^{k-1} ,  ~~~ n=1,\dots,N.
      \label{eq:pararealAlg}
   \end{equation}
\end{enumerate}
Note that the calculation of the fine integrator ${\mathcal F}_H u_{n-1}^{k-1}$ in \eqref{eq:pararealAlg} requires only the initial condition 
$u_{n-1}^{k-1}$, which depends on the previous iteration. Hence, for each $k$,
$\mathcal{F}_t u^{k-1}_{n-1}, $ $0<t \le H$, $n=1,2,\cdots,N$ can be computed in parallel.
The solution computed by the accurate but expensive integrator
is a fixed point. 
Indeed, when the iteration is sufficiently large ($k\ge n$), the solution $u^k_n$ 
become identical to it:
\[
u^k_n = \mathcal{F}_{nH} u_0,~~~n\le k.
\]
In fact, \eqref{eq:pararealAlg} can be regarded as a fixed-point iteration.
In \cite{parareal-Maday10}, it is proved  that under some sufficient conditions of $f$, which we shall recall in Section~\ref{sec:convergence-proof-parareal}, 
\begin{equation}
   | u_n^k - u (nH) | \le C ( H^k + E_f ),
\label{eq:generalErrorEstimate}
\end{equation}
where $E_f $ is the global error in solving the full ODE using the fine 
propagator, and 
$C$ depends on the derivatives of the solutions. 
Eq. \eqref{eq:generalErrorEstimate} assumes a 1st order coarse integrator.
%

In order to identify the source of the difficulty in developing parareal algorithms for highly oscillatory problems, we adapt the
parareal proof of convergence given by Maday in \cite{parareal-Maday10}
for non-singular ODEs.

\subsection{Convergence of parareal}
\label{sec:convergence-proof-parareal}
We consider ODEs of the form \eqref{eq:generalODE}
with initial conditions $u(0)=u_{0}\in D\subset\mathbb{R}^{d}$. 
We are interested in solving \eqref{eq:generalODE} in a fixed time segment $[0,T]$.
The solution is denoted $u(t;u_0),\;t\in [0,T]$.
Let $\Phi$ denote the flow map (propagator) associated with \eqref{eq:generalODE},
\[
   {\Phi}_t x = u(t;x),~~ \forall t>0. 
\]
For sufficiently smooth $f$ we have that
%
   $  \left| {\Phi}_t x -  {\Phi}_t y \right| \le  e^{C t} |x-y|$.
%
In the following $C$ denotes a generic positive constant which may
depend on $T$. 
Since $t \le T$, the prefactor $e^{C t}$ can be bounded by $1 +T^{-1} e^{CT} t$.
This yields a linear stability bound for ${\Phi}_t$,
  \begin{equation}
      \left| {\Phi}_t x -  {\Phi}_t y \right| \le (1+ C t)|x-y| .
\nonumber
\end{equation}

For simplicity, we assume that the coarse integrator ${\mathcal C}_t$ is a one-step 
method with step size $H$ while 
the fine integrator ${\mathcal F}_t$ has step size $h\ll H$.
In addition, we make the following accuracy and stability assumptions on the numerical integrators:
\begin{equation}
       \left| {\mathcal F}_t  x-  {\Phi}_t  x \right| \le C t E_f (1+ |x|),\;\;
       \left| {\mathcal C}_t  x-  {\Phi}_t  x \right| \le C t E_c (1+|x|)  
\label{app:finiteOrder}
\end{equation}
where $E_f$ and $E_c$ denote the global sup error in solving 
\eqref{eq:generalODE} in $[0,T]$
using respectively the fine and coarse integrators in the entire domain of interest $D$.
Note that both $E_f$ and $E_c$ typically depend on $T$. 
In addition,
\begin{equation}
       \left| {\mathcal F}_t x -  {\mathcal F}_t y \right| \le (1+t C)|x-y|, \;\;
       \left| {\mathcal C}_t x -  {\mathcal C}_t y \right| \le (1+t C)|x-y|  
\label{app:stable1a}
\end{equation}
%

%
Let $\delta {\mathcal F}_t = {\Phi}_t - {\mathcal F}_t$ and 
$\delta {\mathcal C}_t = {\Phi}_t - {\mathcal C}_t$,
denote the errors in the fine and coarse propagators, respectively.
Then, by a triangle inequality,
\begin{equation}
       \left| \delta {\mathcal F}_t x - \delta {\mathcal F}_t y \right| \le (1+t C  E_f )|x-y|,\;\;
       \left| \delta{\mathcal C}_t x -  \delta {\mathcal C}_t y \right| \le (1+t C E_c)|x-y|  .
\label{app:stable2}
\end{equation}
We recite the following theorem from \cite{parareal-Maday10}.

\begin{thm}

\label{thm:convergence}
  Let $K \le N/2 = T/2H$. Then, for all $k \le K$,
  \begin{equation}
     \sup_{n=0, \dots, N} | u_n^k - {\mathcal F}_{nH} u_0 | \le C  (E_c)^k .  
\nonumber
  \end{equation}
Consequently, 
\begin{equation}
\sup_{n=0, \dots, N} | u_n^k - u(nH) | \le C  \left[ (E_c)^k + E_f \right] .  
\label{eq:pararealEstEc}
\end{equation}
\end{thm}

\noindent
{\bf Proof:}
   Applying the parareal iterations \eqref{eq:pararealAlg},
\begin{equation}
\begin{aligned}
      u_{n}^{k}- {\mathcal F}_{nH} u_0 
      =& \left[ {\mathcal C}_H u_{n-1}^{k} - {\mathcal C}_H u_{n-1}^{k-1} \right] +  
           \left[ {\mathcal F}_H u_{n-1}^{k-1} - {\mathcal F}_H {\mathcal F}_{(n-1)H} u_0 \right] \\
		=& \left[ {\mathcal C}_H u_{n-1}^{k} - {\mathcal C}_H {\mathcal F}_{(n-1)H} u_0 \right]
		- \left[ \delta {\mathcal C}_H \left({\mathcal F}_{(n-1)H} u_0\right) - \delta {\mathcal C}_H u_{n-1}^{k-1} \right] \\
		& - \left[ \delta {\mathcal F}_H u_{n-1}^{k-1} - \delta {\mathcal F}_H \left({\mathcal F}_{(n-1)H} u_0\right) \right] ,
\end{aligned}
\label{eq:parareal_equalities}
\end{equation}
%
Using assumption \eqref{app:stable2},
%
and denoting \(       \theta_n^k =C (1+CH)^{k-n}  ( E_f + E_c )^{-k} H^{-k} 
       \left| u_{n}^{k} - {\mathcal F}_{nH} u_0 \right|, \)
%
%
%
we have
  $ \theta_n^k \le \theta_{n-1}^k + \theta_{n-1}^{k-1} $.
%
By induction, $\theta_n^k \le C \left(\begin{array}{c} n\\ k \end{array}\right)$.
%
%
Assuming that $E_f < E_c$, 
\begin{equation}
   \left| u_{n}^{k} - {\mathcal F}_{nH} u_0 \right| 
\le C T^K E_c^k 
   = {\mathcal O} (E_c^k).
\label{app:proofend}
\end{equation}
%

\subsection{Parareal and HiOsc ODEs}
\label{sec:notWorking}

We consider HiOsc ODEs given in the singular perturbation form
\begin{equation}
   \dot{u}=\epsilon^{-1}f_{1}(u)+f_{0}(u),
   \label{eq:HiOscODEgeneral}
\end{equation}
with initial condition $u(0)=u_{0}\in D\subset\mathbb{R}^{d}$, where
$D$ is a domain uniformly bounded in $\epsilon$. 
The parameter $0<\epsilon\leq\epsilon_{0}\ll1$
characterizes the separation of time scales -- the fast scale involves
oscillations with frequencies of order $\epsilon^{-1}$ while the
computational time domain is $[0,T]$ with $T$ independent of $\epsilon$.
Throughout the paper we assume that  $f_{1},\; f_{0}$ are sufficiently smooth, and that for each $u_0\in D$, $u(t)$ is  uniformly bounded in $\epsilon$ in the time interval $[0,T]$. 
Furthermore,
we assume that the Jacobian of 
$f_1$ has only purely imaginary eigenvalues in $D$, which are bounded away from $0$
and independent of $\epsilon$.
These settings typically imply that the computational complexity
of direct non-multiscale methods is at least $\mathcal{O}(\epsilon^{-1})$.

To understand some of the challenges in applying the parareal framework
to HiOsc systems, we consider the following simple example 
%
\begin{equation}
  \dot{u} = (\alpha+i \ei )u,~~~u(0)=1 .
\label{eq:simpleSpiral}
\end{equation}
With $\alpha>0$, the trajectory of  $u(t) = e^{(\alpha+i \ei) t}$ is an expanding spiral in the complex plane.
We further assume that the fine integrator is exact, 
${\mathcal F}_t u= e^{(\alpha+i \ei) t} u$.
We first investigate the performance of Algorithm \ref{alg:parareal} using two conventional methods as $\mathcal{C}_t$: Implicit Euler, Explicit Euler, and Trapezoidal Rule.
Table~\ref{table:naiveSpiral} compares the minimal number of parareal iterations, $K$, 
to reach an absolute error  below $1/10$. 
We observe that when conventional methods are implemented as a coarse integrator 
$K$  becomes prohibitively large as $\ep$ gets small.
The increase in $K$ for conventional coarse integrators can be explained by the error estimate \eqref{eq:generalErrorEstimate}.
The difficulty lies in the constant $C$, which grows rapidly with $1/\ep$.
For an order $p$ coarse integrator, the error is proportional to the
$p+1$ time derivative of $f$, which is of order $\mathcal{O} (\ep^{-(p+1)})$.
As a result, the parareal error for HiOsc systems \eqref{eq:pararealEstEc} depends on $\ep$,
\begin{equation}
   | u_n^k - u (nH) | \le C \left\{ E_f +\left[
    \ei \left( \ei H \right)^{p}  \right]^k \right\}.
\label{eq:parareal_error_singular}
\end{equation}
An immediate consequence is that $H$ has to be ${o}(\epsilon)$, 
even when applying A-stable or symplectic methods.
See for example the conclusion in \cite{parareal-MadayLegoll10}.

\begin{table}
\begin{centering}
\begin{tabular}{|c|c|c|c|c|c|c|c|c|c|}
\hline 
Coarse integrator & $\epsilon$= & 0.2 & 0.1 & 0.05 & 0.02 & 0.01 & 0.001\tabularnewline
\hline 
Explicit Euler $(H=\ep/5)$ & $K$  & 7 & 12 & 22 & 52 & 607 & 12200\tabularnewline
Explicit Euler $(H=1/10)$ & $K$  & 34 & 79 & 100 & 100 & 100 & 100 \tabularnewline
\hline 
Implicit Euler $(H=\ep/5)$ & $K$ & 6 & 8 & 13 & 25 & 44 & 351\tabularnewline
Implicit Euler $(H=1/10)$ & $K$ & 18 & 49 & 93 & 100 & 100 & 100\tabularnewline
\hline 
Trapezoidal Rule $(H=\ep/5)$ & $K$ & 1 & 1 & 2 & 3 & 5 & 29\tabularnewline
Trapezoidal Rule $(H=1/10)$ & $K$ & 4 & 18 & 71 & 100 & 100 & 100\tabularnewline
\hline 
The proposed  method $(H=1/10)$ & $K$  & 1 & 1 & 1 & 1 & 1 & 1\tabularnewline
\hline 
\end{tabular}
\par\end{centering}

\caption{\label{table:naiveSpiral}The number of parareal iterations required
to yield an absolute errors of $1/10$ in the expanding spiral example. 
Parameters are $\alpha=1/10$, $T=10$. The maximal number of 
iterations is $T/H$.
}
\end{table}
This simple example reveals the reason why a naive implementation of the
parareal approach may not be effective for integrating HiOsc problems: both
stability and accuracy restrictions require that  the coarse integrator  take steps of order $\ep$.
As a result, the number of coarse steps is $\mathcal{O}(\ei)$ and the method 
may take ${\mathcal O}(\ep^{-1})$ iterations to converge.
For comparison,  we also include in Table~\ref{table:naiveSpiral} the results obtained using the proposed multiscale parareal method.  

\subsection{Layout} 
The layout of the paper is as follows. 
\noindent
Given a conventional parareal method in \eqref{eq:pararealAlg},
\[
\ensuremath{u_{n}^{k}={\mathcal{C}}_{H}u_{n-1}^{k}+{\mathcal{F}}_{H}u_{n-1}^{k-1}-{\mathcal{C}}_{H}u_{n-1}^{k-1}},
\]
Section~\ref{sec:hiosc} presents the main difficulty in using a multiscale method as the coarse integrator in the parareal framework. 
Section~\ref{sec:generalApproach} suggests a general approach for overcoming this difficulty. In particular, two versions of the update using the fine solution will be proposed:
\bigskip \newline
Algorithm~\ref{alg:pararealSlow}: {\it Jacobi style} update for approximating slow variables,
\begin{center}
\includegraphics[width=8.5cm]{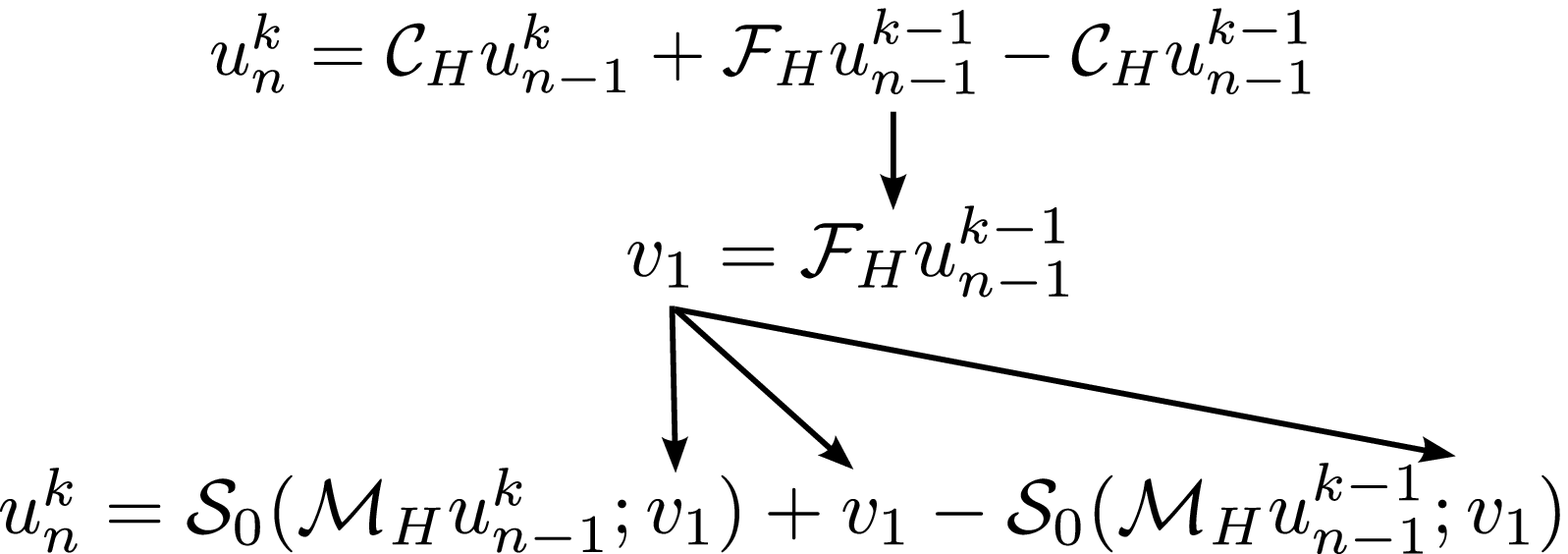}
\par\end{center}
Algorithm~\ref{alg:pararealFull}: {\it Gauss-Seidel style} update for approximating a state variable,

\begin{center}
\includegraphics[width=8.5cm]{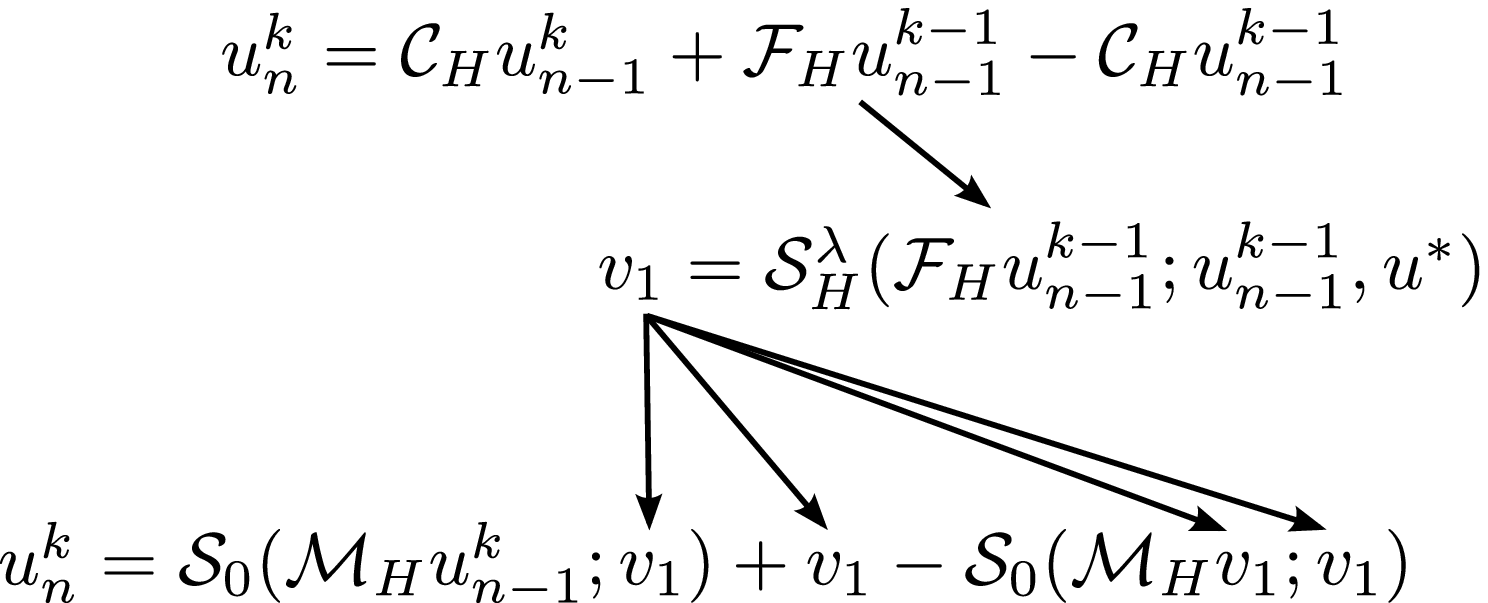}
\par\end{center}
where $u^{*}$ is an approximation of $u((n-1)H)$
from the previous step.
Section~\ref{seq:HiOscParareal} describes the implementations of this strategy:
${\mathcal{S}}_{0}(u_{0};v_{0})$ -- {\it local alignment} in Section~\ref{subsec:local_alignment} and 
$\mathcal{S}_{H}^{\lambda}(u_{1};u_{0},v_{0})$ -- {\it forward alignment} in Section~\ref{subsec:forward_alignment}. 
Section~\ref{sec:reviewPoincare} reviews the Poincar\'e method,
a multiscale numerical method for efficient integration of HiOsc
ODEs presented in \cite{BFHMM2012}.
This method will be used as a coarse solver in the numerical examples presented in
Section~\ref{sec:examples}.
We conclude in Section~\ref{sec:summary}.

\section{Fast oscillations and parareal}
\label{sec:hiosc}
\setcounter{equation}{0}

In order to facilitate the presentation of the main algorithms, we shall first describe the
setting for the underlying multiscale methods.

The literature on efficient numerical integration of problems with separated
time scale is rapidly growing.
For HiOsc ODEs, recent approaches include
envelope methods \cite{Petzold:1997},
FLow AVeraging integratORS \cite{flavors},
Young measure \cite{Artstein-Kevrekidis-Slemrod-Titi07,Artstein-Linshiz-Titi07}
and equation free approaches \cite{Equation-free-review09},
Magnus methods \cite{casas_explicit_2006,iserles_think_2002},
Filon methods \cite{iserles_quadrature_2004,khanamiryan_quadrature_2008},
spectral methods \cite{iserles_numerical_2004,levin_fast_1996},
asymptotic expansions \cite{condon_second-order_2010,iserles_efficient_2005}
and the Heterogeneous Multiscale Methods \cite{HMM-acta-numerica12,E-Engquist-HMM:2003,HMM-review07}.
For a recent review see \cite{engquist_highly_2009}.

Typically, multiscale methods tackle the computational difficulty 
in solving HiOsc ODEs by taking advantage of scale separation, and aim at computing
only the slowly varying properties of the oscillatory solutions.
It requires that enough information about the influence of fast scales
on the slower scale dynamics can be obtained by performing
localized simulations over short times, and thereby better efficiency is  achieved.
The numerical complexity of these methods is therefore
much smaller than direct simulations of the given systems with HiOsc solutions.
For example, \cite{Ariel-Engquist-Tsai07} presents multiscale
algorithms that compute the effective behavior of HiOsc
dynamical systems by using slow variables that are predetermined 
either analytically or  numerically.
More precisely, we define a slow variable
for the system \eqref{eq:HiOscODEgeneral} with solution $u(t;\ep)$ as follows.

\begin{definition}
\label{def:slowvariable}
   A smooth function $a(t,\epsilon)$ is called \textbf{slow} to order $\nu\ge 1$ 
   if $|d^\nu a / d t^\nu|\leq C$ in $t\in[0,T]$ 
   for some constants $C$ and $T$ independent of $\epsilon\in(0,\epsilon_0]$, $\ep_0>0$.  
    A smooth function $\xi(u):D\rightarrow\mathbb{R}$
   is called a \textbf{slow variable} with respect to $u(t)$ if $\xi(t)=\xi(u(t;\ep))$ is slow to order~$1$.
\end{definition}

\noindent
See \cite{Ariel-Engquist-Tsai08,Artstein-Kevrekidis-Slemrod-Titi07,Artstein-Linshiz-Titi07,Fatkullin-Vanden,Gear-Kevrekidis-constraints:2005,Kreiss:1979,Kreiss-Acta-Numerica,Kreiss-Lorenz} for similar definitions and applications
for HiOsc problems.

In this paper, we will work with the following main assumption.
\begin{assumption}\label{assumption:canonical-transform}
There exists a diffeomorphism  $\Psi : u \rightarrow (\xi(u),\phi(u))$,
independent of $\ep$,
separating slow and fast variables such that $(\xi,\phi)$ along the trajectories of \eqref{eq:HiOscODEgeneral}  satisfies an ODE of the form
\begin{equation}
   \begin{cases}
      \dot{\xi}=g_{0}(\xi,\phi), & \xi(0)=\xi(u_0),\\
      \dot{\phi}=\epsilon^{-1}g_{1}(\xi)+g_{2}(\xi,\phi), & \phi(0)=\phi(u_0),
   \end{cases} 
\label{eq:odefastslow} 
\end{equation}
where $\xi\in\mathbb{R}^{d-n}$, $\phi\in\mathbb{R}^{n}$, and $0<\epsilon\le\epsilon_0\ll 1$
is a small parameter. 
We assume that for fixed slow coordinates $\xi$, the fast variable $\phi$ 
is ergodic 
\footnote{By ergodic, we mean that any trajectory of $\phi(t)$
	can get arbitrarily close to any point in the invariant manifold. In particular, this implies the existence of a unique invariant distribution and Birkhoff's ergodic theorem.}
with respect to
an invariant manifold,
which is diffeomorphic to a $n$-dimensional torus, $\mathbb{T}^{n}$.
\end{assumption}
\noindent
Using ergodicity, one can invoke a theory of averaging \cite{Sanders-Verhulst:Averaging}, which implies that the dynamics of slow variables can be approximated ($\mathcal{O} (\ep)$ in the sup norm for $0\leq t \leq T$, $T=\mathcal{O}(1)$) by an averaged equation of the form
\begin{equation}
\begin{aligned}
   & \dot{\bar \xi} = F(\bar \xi),~~~ F({\bar \xi}) = \int g_0(\xi,\phi) d \phi_\xi,\\
   & \bar\xi(0)=\xi(u_0),
   \label{eq:averagedEqn}
\end{aligned}
\end{equation}
where $d \phi_\xi$ denotes the invariant measure for $\phi$ at fixed $\xi$.
For example, perturbed integrable Hamiltonian systems constitute a wide class of dynamical
systems that satisfy this assumption. 
From now on, we shall refer to $\phi$ as the phase of $u$.

The main objective of many multiscale methods is efficient numerical approximations of $\xi(u(t))$ only.
The general strategy of our algorithm is based on 
such multiscale methods for HiOsc ODEs 
that only resolve the macroscopic behavior of a system as
specified by the slow variables 
\cite{Iterated_averaging,tutorial,Ariel-Engquist-Tsai07,Ariel-Engquist-Tsai08,AET-reversible,AET-threeScale,Engquist-Tsai-HMM-ODE:2003,Sharp05}.
In this respect, the algorithms listed above are different from other multiscale methods that resolve all scales of the dynamics, 
for example, multi-level methods or high-order asymptotic expansions \cite{casas_explicit_2006,condon_efficient_2012,condon_second-order_2010,Kreiss:1979,Kreiss-Acta-Numerica}.


It is possible to design a parareal algorithm for computing {\it only the averaged 
slow variables}
using multiscale integrators as both the coarse and fine integrators.
Such an approach is essentially a parareal scheme for the averaged equation.
However, this is not the point of this paper --- 
here we are interested in the possibility of creating a parareal
algorithm that computes {\it all state variables}, including the fast phase information.

We consider the problem of using a multiscale integrator in the coarse integration, and 
provide the stability of the corresponding coupling of multiscale-fine integrators under the parareal framework. 
Since the error bound stated in \eqref{eq:parareal_error_singular} still formally applies in this case,
one cannot expect convergence of $u(t)$ unless some additional improvement is made to the chosen existing multiscale scheme.

Consider again the simple expanding spiral \eqref{eq:simpleSpiral} with $\alpha=1$.
It is easily verified that $|u(t)|=e^t$ is a slow variable.
For convenience of the discussion, we assume
that the fine/microscopic solver is exact, i.e. ${\mathcal F}_t u=e^{(1+i\ei)t}u$,
and that the coarse/macroscopic solver is \emph{exact in the slow variables,
}i.e. any function of $|u|$ is computed without error \emph{but
the phase of $u$ may be wrong.} 
We write the
macroscopic solution as $\mathcal{C}_{t}U=e^{t}e^{i(t\epsilon^{-1}+\theta_{t})}U$,
where $\theta_{t}\in[0,2\pi)$ denotes the error in the phase that
is produced by the macroscopic solver. 
Applying Algorithm~\ref{alg:parareal} we obtain 
\begin{equation} 
u_{3}^{(1)}  =u(3H)\left(1+\mathcal{O}(\theta_{H}^{2})\right), ~~
u_{3}^{(2)}  =u(3H)\left(1+\mathcal{O}(\theta_{H}^{3})\right).
\end{equation} 
This simple exercise shows that \emph{the naive
iterations improve the accuracy of the macroscopic solution if $\theta_{H}$
is small}, \emph{and that the iterations diverge if $\theta_{H}$
is not sufficiently small.} However, in a typical HiOsc, $\theta_{H}$ is
not necessarily small. In general, $\theta_{H}$ can be any value in $[0,2\pi]$
and is not necessarily small.

In the following sections, we show that by aligning the phase of the coarse and fine solvers, 
it is possible to design parareal algorithms that use multiscale coarse integrators.

\section{Multiscale parareal}
\label{sec:generalApproach}
\setcounter{equation}{0}

In this section, we introduce the main contribution of this paper -- accurate and convergent parareal algorithms that use multiscale methods as coarse integrators.
Two parareal schemes are presented. The first focuses on approximating only the slow variable, 
while the second achieves sup-norm convergence in the state variable, $u\in\mathbb{R}^d$. 
Both methods are based on a phase alignment  strategy, which can be applied 
if, for fixed slow variables, the phase is ergodic with respect to a circle. 
Accordingly, we assume that the slow coordinate $\xi=(\xi_1, \cdots, \xi_{d-1})$ is a vector of $d-1$ functionally independent slow variables.

\subsection{Multiscale coarse integrator}
\label{sec:onlySlow}

For the remainder of this paper, we shall assume that  the coarse propagator
is a multiscale method that only approximates the slow variables.
In order to emphasize this point, the multiscale coarse integrator will be denoted 
${\mathcal M}_t$ in place of ${\mathcal C}_t$.
Similar to assumption~\eqref{app:stable2}, we shall assume that 
\begin{equation}
\begin{aligned}
	\left|\delta {\mathcal M}_t x - \delta {\mathcal M}_t y \right| & \le (1+t C H) |\xi(x)-\xi(y)|.
\end{aligned}
\label{eq:doesNotHold}
\end{equation}
The parareal proof of convergence as given in 
Section~\ref{sec:convergence-proof-parareal} hinges on the stability assumption \eqref{app:stable2},
which does not directly involve the exact solution.
As a result, as long as \eqref{app:stable2} holds, the parareal iterations
will converge, although not necessarily to the exact solution.
However, with a multiscale coarse integrator, \eqref{eq:doesNotHold} implies that stability only in the slow variables is guaranteed.
Accordingly, we propose to  
modify the coarse multiscale integrator by fixing the fast variable 
(a fast phase in the case of HiOsc problems).
In terms of slow-fast coordinates, the multiscale integrator will be stable in the slow coordinates
due to \eqref{eq:doesNotHold} while stability in the fast variable will be
enforced by aligning trajectories with respect to a common reference phase.
In order to achieve this, we assume that one can devise the following local alignment
algorithm.

\vspace{0.3cm}

\noindent 
\centerline{
\framebox{\begin{minipage}[c]{1\columnwidth}%
{\it \textbf{Local alignment: }}\medskip \\
Given $u_{0}$ and $v_{0}$ such that $\xi(u_{0})=\xi(v_{0}) + \Delta \xi$. \\
Let $w_0 = \Psi^{-1} (\xi(u_0), \phi(v_0))$ be the point that has the same slow 
coordinates as $u_0$ and the same phase as $v_0$. \medskip \\
Find a point $\tilde{w}_0$ such that $\left| \tilde{w}_0 - w_0 \right| = {\mathcal O} (\Delta \xi)$.
\end{minipage}}
}
\vspace{0.3cm}

\noindent 
In other words, the local alignment procedure replaces $u_0$ by a new point $\tilde{w}_0$ that has
the same (to order $\Delta \xi$) slow coordinates, i.e. $\xi$ values, as $u_0$, and approximately the same phase as $v_0$. 
A trivial solution to the local alignment problem is to set $\tilde w_0:= v_0$. 
However, this is not an adequate strategy that can be used in the next steps of development 
of our multiscale parareal algorithm. 
\begin{notation}
We denote such a local alignment procedure as $\tilde{w}_0={\mathcal S}_0 (u_0;v_0)$.
\end{notation}

Given a local alignment algorithm ${\mathcal S}_0$, we propose the following 
modified parareal scheme.
\begin{algorithm}
\label{alg:pararealSlow}
\end{algorithm}
\begin{enumerate}
\item Initialization: (Construct the zero'th iteration approximation)
   \begin{align}
       u_0^0 =& u_0 ~~ {\rm and} ~~ u^0_{n} = {\mathcal M}_H u^0_{n-1}, ~~~ n=1,\dots,N.
\nonumber
   \end{align}
\item Iterations: $k=1 \dots K$
   \begin{enumerate}
      \item Parallel fine integrations for $n=k,\dots,N$,
           \begin{equation}
              u_{F,n}^k = {\mathcal F}_H u^{k-1}_{n-1} . 
\nonumber
           \end{equation}
   \item Parareal correction: For $n=k,\dots,N$,
   \begin{equation}
       u_0^k = u_0 ~~ {\rm and} ~~ 
      u_{n}^{k} = {\mathcal S}_0 ({\mathcal M}_H u_{n-1}^{k} ; u_{F,n}^{k}) + u_{F,n}^{k} -    
       {\mathcal S}_0 ({\mathcal M}_H u_{n-1}^{k-1} ; u_{F,n}^{k}).
\label{eq:pararealAlg3}
   \end{equation}
\end{enumerate}
\end{enumerate}
In each iteration we first calculate all fine scale integrations.
Then, the results of the multiscale integrators are aligned with the fine scale ones.
In the following, we prove that using Algorithm~\ref{alg:pararealSlow},
all slow variables converge to their limiting value given by the fine scale approximation. 
We consider a 1st order multiscale integrator with local phase alignment. 

\begin{thm}
\label{thm:slow_convergence}
  Let $K \le N/2 = T/2H$. Then, for all $k \le K$,
\begin{equation}
\sup_{n=0,\dots,N} \left| \xi (u_{n}^{k} ) - \xi \left( {\mathcal F}_{nH} u_0 \right) \right| \le  CH^k .
\nonumber
\end{equation}
\end{thm}

\noindent
{\bf Proof:} We recall the assumption that there exists a diffeomorphism $\Psi:u \to (\xi(u),\phi(u))$
such that $\xi\circ u(t)$ are slow while $\phi\circ u(t)$ are fast. 
The variables $(\xi,\phi)$ are only used in the analysis but {\em not} in the numerical algorithm.

The main difference with the general analysis described in Section~\ref{sec:convergence-proof-parareal} 
is that the bound \eqref{app:stable2} is not valid if a multiscale coarse integrator is used.
Instead, denoting by
$\delta {\mathcal S}_0 ({\mathcal M}_t u_1; u^*)={\mathcal S}_0 ({\Phi}_t u_1; u^*) -  {\mathcal S}_0 ({\mathcal M}_t u_1; u^*)$,
we have
\begin{equation}
   \Psi\circ\delta {\mathcal S}_0({\mathcal M}_t u_1; u^*) - 
       \Psi\circ\delta {\mathcal S}_0 ({\mathcal M}_t u_2 ; u^*)
       = \left(\delta\xi, \delta\phi\right), 
  \label{app:assumptions1}
\end{equation}
such that
%
    $ | \delta\xi |  \le (1+C t) |\xi(u_1)-\xi(u_2)| $
but $| \delta \phi |={\mathcal O} (\ep)$ which is the accuracy of
local alignment. 
Comparing with conventional methods as a coarse integrator and the related estimate
\eqref{eq:parareal_error_singular}, 
the slow part is \emph{controlled by the local phase alignment} in Algorithm~\ref{alg:pararealSlow}
just like in the non-singular case, while the rapidly changing phase is
incorrect but does not affect the accuracy of the slow variables.

The slow variables of $u^k_n$  in \eqref{eq:pararealAlg3} are
\begin{equation}
\xi(u^k_n) = \xi ({\mathcal S}_0 ({\mathcal M}_H u_{n-1}^{k} ; u_{F,n}^{k})) + 
\xi (u_{F,n}^{k}) -\xi({\mathcal S}_0 ({\mathcal M}_H u_{n-1}^{k-1} ; u_{F,n}^{k})),
\nonumber
\end{equation}
which is valid with the local alignment.
We may thus think of the multiscale integrator combined with the local alignment as a coarse integrator
with first order accuracy for the slow variables.
For shorthand, we denote by ${\mathcal M}_H u_{n-1}^{k}$ the combined ${\mathcal S}_0 ({\mathcal M}_H u_{n-1}^{k} ; u_{F,n}^{k})$.
The error in the slow variables is evaluated similarly to \eqref{eq:parareal_equalities},
\begin{equation}
\begin{aligned}
 &\xi(u_n^k) - \xi({\mathcal F}_{nH} u_0) \\
     &= \left[ \xi ({\mathcal M}_{H} u_{n-1}^k) 
           + \xi ({\mathcal F}_{H} u_{n-1}^{k-1})
           - \xi ({\mathcal M}_{H} u_{n-1}^{k-1})
           - \xi ({\mathcal F}_{H} {\mathcal F}_{(n-1) H} u_0)
           \right] \\
    &= \left[ \xi ({\mathcal M}_{H} u_{n-1}^k) 
                          - \xi ({\mathcal M}_{H} {\mathcal F}_{(n-1) H} u_0)
           \right] 
           + \left[ 
                          \xi (\delta {\mathcal M}_{H} {\mathcal F}_{(n-1) H} u_0)
                          - \xi (\delta {\mathcal M}_{H} u_{n-1}^{k-1})
           \right] \\
           &~~+ \left[ 
                          \xi (\delta {\mathcal F}_{H} u_{n-1}^{k-1})
                          - \xi (\delta {\mathcal F}_{H} {\mathcal F}_{(n-1) H} u_0)
           \right] 
\end{aligned}
\nonumber
\end{equation}
Using \eqref{app:assumptions1}, 
for every slow variable $\xi$, we have that 
\begin{equation}
\begin{aligned}
      & \left| \xi ( u_{n}^{k} ) - \xi ( {\mathcal F}_{nH} u_0 ) \right| \le \\ 
      & (1 + C H) \left| \xi ( u_{n-1}^{k} ) - \xi ({\mathcal F}_{(n-1)H} u_0 ) \right|
      + C (E_f + H ) H \left| \xi (u_{n-1}^{k-1} ) - \xi ( {\mathcal F}_{(n-1)H} u_0 ) \right|
\nonumber
\end{aligned}
\end{equation}
Denoting
%
     $ \theta_n^k = (1+CH)^{k-n}  ( E_f + H )^{-k} H^{-k} 
      \left| \xi ( u_{n}^{k} ) - \xi \left( {\mathcal F}_{nH} u_0 \right) \right| $
%
and following the same procedure as in \eqref{app:proofend}, we have for the slow variable,
\begin{equation}
\sup_{n=0,\dots,N} \left| \xi (u_{n}^{k} ) - \xi \left( {\mathcal F}_{nH} u_0 \right) \right| \le
   C  ( E_f +  CH)^k \le    CH^k .
\nonumber
\end{equation}

\subsection{Phase continuity in the coarse and fine scale simulations}
\label{eq:phaseCont}

We next consider convergence of the parareal approximation to the exact solutions.
The main idea is to enforce consistency in the fine scale solutions between neighboring coarse time intervals. 
We may rephrase this problem as the following.

\vspace{0.3cm}
\noindent 
\centerline{
\framebox{\begin{minipage}[c]{1\columnwidth}%
{\it \textbf{Forward alignment of step size $H$:}}\medskip\\
Given $u_0$, $v_0$, and $u_1={\mathcal F}_H u_{0}$  
such that $\xi(u_0)-\xi(v_0) = {\mathcal O} (\ep)$. \\
Let $w_0 = \Psi^{-1} (\xi(u_0), \phi(v_0))$ and  $w_1 = {\mathcal F}_H w_0$.\medskip\\
Find a point $\tilde{w}_1$ such that $\xi(\tilde{w}_1)=\xi(w_1) + {\mathcal O} (\ep)$
and $\phi  (\tilde{w}_1)= \phi(w_1) + \mathcal{O} (H^2)$.
\end{minipage}}
}
\vspace{0.3cm}

\begin{figure}[tbh]
   \begin{centering}
   \includegraphics[width=14cm]{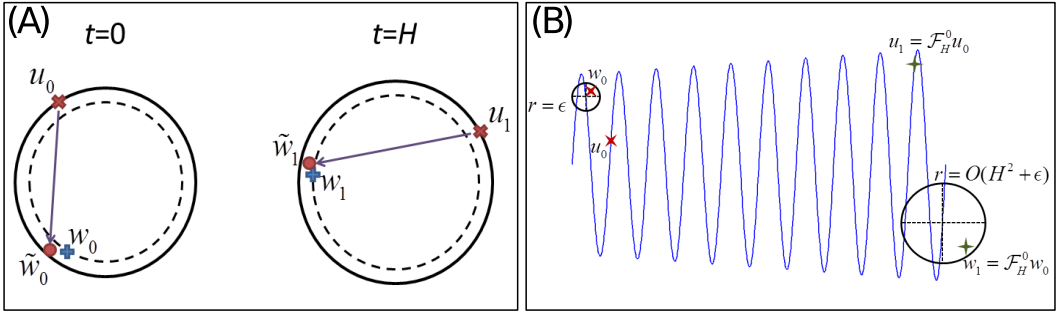}
   \par\end{centering}
   \caption{
   \label{phaseCorrections} 
         Local and forward alignments. 
         (A) At $t=0$, given two points $u_0$ and $v_0$, we wish to approximate
         $w_0 = {\mathcal S}_0 (u_0 ; v_0) $ -- a point that has the same slow variables as $u_0$ and the same phase as $v_0$.
         At $t=H$, we approximate the point $w_1={\mathcal F}_H w_0$
         (B) At $t=0$, a small ${\mathcal O} (\ep)$ of $\tilde{w}_0$ yields
         at $t=H$ a larger ${\mathcal O} (H^2 + \ep)$ error.
	  The center of the circles represent  $\tilde{w}_0$ and  $\tilde{w}_1$.
         }
\end{figure}

\noindent
In the problem of forward alignment, if $w_0$ is a point with the same slow variable as $u_0$ and phase as $v_0$,
then a forward alignment procedure constructs an order $H^2$ approximation of 
$w_1={\mathcal F}_H w_{0}$, the right end point of a coarse interval.
See Figure~\ref{phaseCorrections}A for a schematic sketch.
\begin{notation}
We denote such a forward alignment procedure as
$\tilde{w}_1={\mathcal S}_H^\lambda (u_1;u_0,v_0)$,
where $\lambda$ are precomputed parameters to be used in the alignment.
\end{notation}

The forward alignment procedure can be trivially accomplished simply by setting $\tilde w_1$ to 
$v_1{\mathcal F}_H v_0$ or $w_1$.
However,  this would require the additional computation of $v_1$ from $v_0$, and so this trivial "fix" has a computational cost of sequentially solving the entire system with the fine integrator.
In practice, for the purpose of parallel in time computations, 
one needs to do so with a computational
cost that is lower than running the fine scale solver sequentially.
Hence, we need to estimate the solution to the given ODE with the given initial condition $v_0$ by certain simple
operations performed on the fine scale solutions already computed in parallel.
In the following section, we shall describe a forward alignment algorithm 
for the special case of HiOsc ODEs, in which, for fixed slow variables, 
the fast phase is periodic.
The method applies only a local exploration by means of minimal additional fine scale computation of the solution around $u_0$ and $u_1$.
In particular, its efficiency is independent of $\epsilon$.

To summarize, we present the complete multiscale-parareal algorithm.
Recall that $\mathcal{M}_H$ is a multiscale method that only approximates the slow variables.
For the fast phase,
using local and forward alignments, we propagate the needed phase adjustments sequentially 
along with the parareal correction.

\medskip
\noindent 
\centerline{
\framebox{\begin{minipage}[c]{1\columnwidth}%
\begin{algorithm}
\label{alg:pararealFull}
Full multiscale-parareal algorithm.\end{algorithm}
\begin{enumerate}
\item Initialization: Construct the zero'th iteration approximation:
\begin{equation}
u_{0}^{0}=u_{0}~~{\rm and}~~u_{n}^{0}={\mathcal{M}}_{H}u_{n-1}^{0},~~~n=1,\dots,N.
\nonumber
\end{equation}
\item Iterations: $k=1...K$
\begin{enumerate}
\item Parallel fine integrations for $n=k,\dots,N$:
\begin{equation}
u_{F,n}^{k-1}:=\mathcal{F}_{H}u_{n-1}^{k-1}.
\nonumber
\end{equation}
\item Header: For $n=0,\dots,k-1$, set $u_{n}^{k}=u_{F,n}^{k-1}$.
\item Parareal step: 
Set the initial reference point $u^{*}=u_{k-1}^{k-1}$ and for $n=k,\dots,N$, 
\begin{enumerate}
\item Locally align the previous $u_{n-1}^{k-1}$ with the current reference point $u^{*}$:
\[
\tilde{u}_{n-1}^{k-1}={\mathcal{S}}_{0}(u_{n-1}^{k-1};u^{*}).
\]
\item Align forward to the end of the coarse segment:
\[
{\tilde{u}_{F,n}^{k-1}={\mathcal{S}}^\lambda_{H}(u_{F,n}^{k-1};u_{n-1}^{k-1},u^{*})}.
\]
\item Corrector: 
\begin{equation}
u_{n}^{k}={\mathcal{S}}_{0}({\mathcal{M}}_{H}u_{n-1}^{k};\tilde{u}_{F,n}^{k-1})+\tilde{u}_{F,n}^{k-1}-{\mathcal{S}}_{0}({\mathcal{M}}_{H}\tilde{u}_{n-1}^{k-1};\tilde{u}_{F,n}^{k-1}).
\nonumber
\end{equation}
\item Update the reference point $u^{*}={u}_{n}^{k}$ and repeat.
\end{enumerate}
\end{enumerate}
\end{enumerate}
\end{minipage} 
}
}
\medskip

\noindent
Local and forward alignment steps (Step 2(c)i and ii, respectively) create a point $\tilde{u}_{F,n}^{k-1}$ at the end of each coarse segment
according to which
all points in the current corrector iteration can be aligned.
Since the error in each forward alignment is of order $H^2$, 
we find that the overall phase is continuous up to a global ${\mathcal O} (H)$ error.
We conclude that following phase alignments, 
the aligned coarse multiscale method provides an globally ${\mathcal O} (H)$ approximation
of both slow and fast variables, i.e., it approximates
the solution in the sup norm.

Example~\ref{sec:spiral} demonstrates the effectiveness of the method 
in a more complicated expanding spiral with a slowly changing frequency. 
Before describing numerical methods for local and forward alignments of HiOsc ODEs, 
we address the convergence of the algorithm.

\subsection{Convergence of Algorithm~\ref{alg:pararealFull}}

Convergence of Algorithm~\ref{alg:pararealFull} in the state variable is obtained in two steps.
First, following Section~\ref{sec:onlySlow} and Theorem~\ref{thm:slow_convergence},
all slow variables converge to their values obtained by the fine scale integrators,
\begin{equation}
  \sup_{n=0,\dots,N}  | \xi(u_n^k) - \xi ( u (nH)) | \le C ( H^k + E_f ),
\nonumber
\end{equation}
where $C$ is a constant that is independent of $\ep$.
In particular, if $E_f={\mathcal O} (\ep)$, then, 
following ${\mathcal O}(\log(\ep))$ iterations, 
the error in the slow variables is of order $\ep$.
As a result, after a few (typically one or two) iterations, the assumptions underlying 
forward alignment, that the error in the slow variables is of order $\ep$ holds (more precisely, $\xi(u_0)-\xi(v_0)={\mathcal O} (\ep)$). 
We may thus think of the adopted multiscale combined with the local/forward alignment algorithm as a coarse integrator
with first order accuracy for all state variables.
Hence, the conventional parareal proof of convergence described in 
Section~\ref{sec:convergence-proof-parareal} holds.
More precisely, suppose that, following the phase alignment, the set $\{ w_0,\dots,w_{N} \}$
is computed and updated in every iteration.
The following estimates hold for $j=0,\dots,N$,
\begin{equation}
\begin{aligned}
   \left| \xi(w_j) - \xi(u(jH)) \right| &= {\mathcal O} (E_m), \\
   \left| \phi(w_j) - \phi(u(jH)) \right| &= {\mathcal O} (H), \\
   \left| w_j - u(jH) \right| &= {\mathcal O} (H)+{\mathcal O} (E_m),
\end{aligned}
\nonumber
\end{equation}
where $E_m$ is the error of the aligned multiscale method in approximating 
the slow variables.

%
Assume further the stability properties for the fine and aligned-multiscale coarse propagators
\eqref{app:stable2}.
%
Therefore, after one parareal iteration, 
%
\begin{equation}
   \sup_{n=0,\dots,N} \left| u_n^1 -  u(nH) \right| \le C ( \ei E_m + E_f ).
\nonumber
\end{equation}
After $k$ iterations, and assuming a first order multiscale coarse integrator, 
$E_m = {\mathcal O}(H)$,
\begin{equation}
   \sup_{n=0,\dots,N} \left| \xi(u_n^k )-  \xi(u(nH)) \right| \le C ( H^k + E_f ) ,
\nonumber
\end{equation}
and
\begin{equation}
   \sup_{n=0,\dots,N} \left| u_n^k -  u(nH) \right| \le C ( \ei H^k + E_f ).
\nonumber
\end{equation}
%
{The accuracy of slow variables is improved to by a factor of $H$ per parallel iteration
(compare with the diverging factor of $(\ei H)^k$ in
\eqref{eq:parareal_error_singular})}.
\begin{rem}
In \cite{parareal-Legoll13}, Legoll et al propose a multiscale parareal algorithm
for stiff ODEs in which the fast dynamics is dissipative, i.e., trajectories quickly converge
to a lower dimensional manifolds.
Unlike the HiOsc case, a naive application of the parareal methodology to 
stiff dissipative systems converges. 
However, it is not very efficient and suffers from similar difficulties as discussed earlier.
To circumvent these difficulties, Legoll et al \cite{parareal-Legoll13} suggest a correction step
that allows a consistent approximation of the fast-slow dynamics with parareal.
This work assumes that the system is split into slow and fast variables, or alternatively, that a change of variables that splits the system into such coordinates is given explicitly. 
Essentially, the idea of \cite{parareal-Legoll13} is to set the fast component of the multiscale solver with that obtained from the fine one.
This step may be viewed as a simple alignment method.
Indeed, if the $(\xi,\phi)$ coordinates are know, then the same approach can also be applied
to the HiOsc case.
In contrast, the method presented in the following section is seamless in the 
sense that it does not require 
knowing the slow nor the fast variables.
\end{rem}

\section{Phase alignment strategies}
\label{seq:HiOscParareal}
\setcounter{equation}{0}

In this section, we describe a numerical method for both local and forward alignments
as defined in the previous section for the special case of HiOsc ODEs
in which, for fixed slow variables $\xi$, the dynamics of the
fast phase $\phi$ is periodic.

In the Algorithm~\ref{alg:pararealFull},
$v_0$ and $u_0$ will correspond to $u^k_n$ and 
$u^{k-1}_{n}$, 
the solutions computed at the current and the previous iterations, respectively.
The assumption is that $v_0$ is the more accurate approximation of the solution at the time $t=nH$, particularly in the phase variable.
The goal is that from the available information, $u_0$, $v_0$, and 
$u_1:=\mathcal{F}_H u_0$, we estimate  $v_1:=\mathcal{F}_H v_0$ at $t=(n+1)H$ in order to make correction in the phase of $u_1$.
We also emphasize that in the subsequent time steps, $\mathcal{F}_H u_0$ is always available because of the prior parallel fine integrations.
Now $w_0$, as defined in Section~\ref{eq:phaseCont},  is a point on the same slow coordinates as $u_0$ 
but has the same phase as $v_0$. Consequently $w_1:=\mathcal{F}_H w_0$ is a good estimate of $v_1$.
In this section, we propose a strategy that move $u_0$ to $w_0$,  
and  $u_1 $ to a state that is within $\mathcal{O}(\epsilon)$ to $w_1$, \emph{without computing $\mathcal{F}_H w_0$ or $\mathcal{F}_H v_0$.}
Our goal is to describe a method that finds a point $\tilde{w}_0$ such that 
$|w_0 - \tilde{w}_0|={\mathcal O} (\ep)$ (local alignment)
and a second point $\tilde{w}_1$ such that  $|w_1-\tilde{w}_1|={\mathcal O} (H^2+\ep)$
(forward alignment).

In addition to Assumption~\ref{assumption:canonical-transform}, we assume the following,

\begin{assumption}
\label{assumption:canonical-transform1}
The fast variable $\phi \in \R$ and $g_0(\xi,\phi)$ is 1-periodic in $\phi$. 
\end{assumption}

For fixed $\xi$, the time derivative of $\phi$ may depend on the slow variables, i.e., the periodicity
in time of $g_0 (\xi,\phi(t))$ is of order $\ep$ and depends on $\xi$. 
Accordingly, it is denoted $\ep \tau(\xi)$, where $\tau$ is a smooth, slow function.
Note that this does not mean that the oscillation in the original state variables
are linear because the transformation $\Psi$ is in general nonlinear. 

\subsection  {At time $t$ move $u_0$ closer to $w_0$} 
\label{subsec:local_alignment}

Assuming that $\xi (u_0) - \xi (v_0) = {\mathcal O} (\ep)$, we may use $v_0$
instead of $w_0$, which is not known.
Denote
 \begin{equation*}
    J (t;u_0,v_0) = \left| {\mathcal F}_t u_0 - v_0 \right|^2 .
\end{equation*}
We look for the local minima of $J (t;u_0,v_0)$ closest to $t=0$ (by the periodicity assumption, such local minima exist)
\begin{equation*}
    0 = J' (t;u_0,v_0)(t) = 2 \left( {\mathcal F}_t u_0 - v_0 \right) \cdot \frac{d}{dt} 
    \left( {\mathcal F}_t u_0 - v_0 \right) = 
    2 \left( {\mathcal F}_t u_0 - v_0 \right) \cdot 
    \left( \frac{\partial \Psi^{-1}}{\partial \xi} \dot{\xi} + \frac{\partial \Psi^{-1}}{\partial \phi} \dot{\phi} \right) .
\end{equation*}
To leading order in $\ep$, we have 
\begin{equation}\label{eq:local_min}
({\mathcal F}_t u_0 - v_0) \cdot (\partial \Psi^{-1}/\partial \phi) = {\mathcal O} (\ep).
\end{equation}
In other words, the phase  of ${\mathcal F}_t u_0$ is close to that of $v_0$,
$\phi ( {\mathcal F}_t u_0 ) = \psi_0 + {\mathcal O} (\ep)$
and therefore also to the phase of $w_0$.
We denote the ``first'' two local minima
\begin{equation*}
   - \ep \tau_0 + {\mathcal O} (\ep^2)< t_0^- < 0 < t_0^+ < \ep \tau_0 + {\mathcal O} (\ep^2) ,
\end{equation*}
where $\tau_0=\tau (\xi_0)$.
Consider
 \[  
   \tilde{w}_0^\pm = {\mathcal F}_{t_0^\pm} u_0,
\]
%
and the convex combination using these points,
%
\begin{equation*}
   \tilde{w}_0 = {\mathcal S}  u_0 = \lambda_+ \tilde w _0 ^+ + \lambda_- \tilde w_0^-  ,
\end{equation*}
with weights $\lambda=(\lambda_+ , \lambda_- )$ independent of $\ep$. 
Thus, equation \eqref{eq:local_min} implies that,
for any linear combination such that 
$\lambda_+ + \lambda_- = 1$, $| \tilde{w}_0 - v_0|=\mathcal{O}(\epsilon)$, i.e.,
$\tilde w_0$ defined above is a valid choice in the local alignment procedure.
We define

\medskip
\noindent 
\centerline{
\framebox{\begin{minipage}[c]{1\columnwidth}%
{\it \textbf{Local alignment: }}
\begin{equation}
{\mathcal S}_0 (u_0;v_0) = \lambda_+ {\mathcal F}_{t_0^+} u_0 + \lambda_- {\mathcal F}_{t_0^-} u_0 .
\nonumber
\end{equation}
\end{minipage}}
}
\medskip

\noindent
In the numerical implementation of the local alignment ${\mathcal{S}}_{0}$,
we use a simple algorithm which solves
the $l_2$ minimization of $J(t)$ with a small ${\mathcal O}(\epsilon)$ step size and adaptively increasing the search domain
until an appropriate minimum is found. Denoting the time where the minimum is attained by $t^{*}$, 
we improve the accuracy by a quadratic interpolation through the neighboring points of $t^{*}$. This algorithm achieves ${\mathcal O}(\epsilon)$ accuracy with an efficiency that is independent of $\epsilon$. 
See Appendix \ref{appen:search_algorithm} for details and \cite{Book:Nocedal} for further references on other efficient minimization techniques.



\subsection{At time $t+H$ move $u_1$ closer to $w_1$}
\label{subsec:forward_alignment}

We would like to do the same at $t=H$, i.e., move $u_1={\mathcal F}_H u_0$ to 
$w_1={\mathcal F}_H w_0$.
The main difficulty is that we cannot expect that the solution has oscillations of constant periodicities.
We denote $\tau_1 = \tau (\xi_1) = \tau_0 +(\xi_1 - \xi_0) (\partial \tau/ \partial \xi ) + {\mathcal O} (H^2) $. 
In analogy to the procedure at $t=0$, we find the ``first'' two minimizers of
 \begin{equation*}
    J (t;u_1,w_1) = \left| {\mathcal F}_t u_1 - w_1 \right|^2 ,
\end{equation*}
such that
\begin{equation*}
    - \ep \tau_1 + {\mathcal O} (\ep^2) < t_1^- < 0 < t_1^+ < \ep \tau_1 + {\mathcal O} (\ep^2) .
\end{equation*}
%
%
Let 
\begin{equation*}
   \tilde{w}_1 = \lambda_+  \tilde w_1^+ + \lambda_- \tilde w_1^-,~\textrm{with}~\tilde w _1^\pm = {\mathcal F}_{t_1^\pm} u_1 = w_1 + {\mathcal O} (\ep).
\end{equation*}
Then, for any  constants $\lambda_+,\lambda_-$  we have that $|w_1 - \tilde{w}_1 |={\mathcal O}(\ep)$.
\emph{The problem is that we do not know $w_1$ and therefore cannot find $t_1^\pm$}. 
One option is to use $t_0^\pm$ instead and choose weights $\lambda_\pm$ that minimize the error.
This requires us to relate the $t_0$'s and $t_1$'s. 

\begin{figure}[tbh]
   \begin{centering}
\includegraphics[width=13cm]{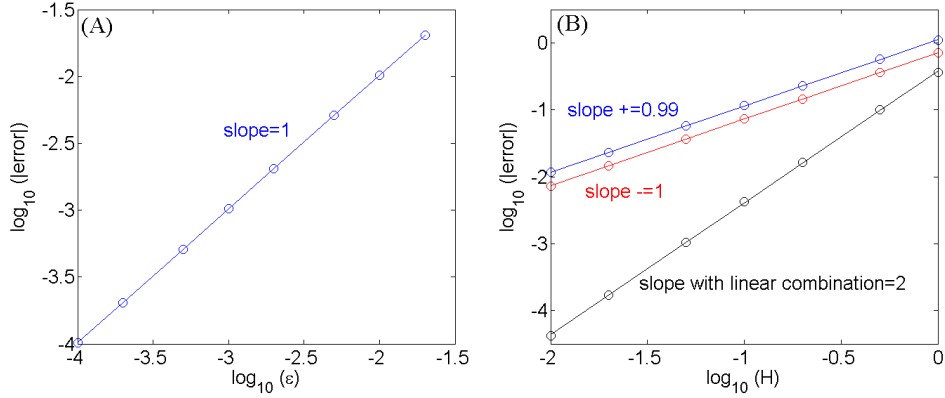}  
   \par\end{centering}
\caption{
  The error in correcting the phase at the end of one coarse segment for the expanding spiral \eqref{eq:simpleSpiral}.
  (A) In local alignment, the phase at the end of a coarse segment $u_0$
  is aligned with $v_0$ with an ${\mathcal O} (\ep)$ error.
  (B) Blue: forward correction with $s_+$, red: backward correction with $s_-$ and
  black: a linear combination of shifts using the forward alignment algorithm defined in 
  Section~\ref{seq:HiOscParareal}. With the proposed linear combination, the error is of order $H^2$.
}
\label{fig:Stilde}
\end{figure}

Denote
\begin{equation*}
\begin{aligned}
   \Psi u_0 = \left( \xi_0, \phi_0 \right), ~~~ & ~~~ 
        \Psi u_1 = \Psi {\mathcal F}_H u_0 =
        \left( \xi_1, \phi_1 \right), \\
   \Psi v_0 = \left(\eta_0, \psi_0 \right), ~~~ & \\
   \Psi w_0 = \left( \xi_0, \psi_0  \right),  ~~~ &
        ~~~   \Psi w_1 = \Psi{\mathcal F}_H w_0 = \left( \xi_1, \psi_1\right).
\end{aligned}
\label{eq:slowFastCoordinates}
\end{equation*}
Without loss of generality, we assume that $\psi_0 > \phi_0$ and $|\psi_0-\phi_0|<1$.
Similarly, assume $\psi_1 > \phi_1$ and $|\psi_1-\phi_1| < 1$.
Then, to leading order in $\ep$,
\begin{equation*}
\begin{aligned}
   t_0^+ = (\psi_0 - \phi_0) \ep\tau_0 ,   ~~~~ & ~~~~  t_1^+ = (\psi_1 - \phi_1) \ep\tau_1, \\
   t_0^- = - (1-\psi_0 + \phi_0) \ep\tau_0 ,   ~~~~ & ~~~~  t_1^- = -(1 - \psi_1 + \phi_1) \ep\tau_1. 
\end{aligned}
\end{equation*}
Next, denote the solution of $(\xi,\phi)$ with initial condition $(\xi_0,\phi_0)$ as
 $\xi (t;\xi_0,\phi_0)$ and $\phi (t;\xi_0,\phi_0)$, i.e., $\left(  \xi (t;\xi_0,\phi_0), \phi (t;\xi_0,\phi_0) \right) =      \Psi {\mathcal F}_t u_0 $.
Using the averaging principle \eqref{eq:averagedEqn}, we can write
\begin{equation*}
   \xi (t;\xi_0,\phi_0) = \bar\xi (t) + \ep \gamma (t/\ep,\xi)  + {\mathcal O} (\ep^2) ,
\end{equation*}
where $\bar\xi (t )$ is a slow function that does not depend on the phase 
and $\gamma(s,\xi)$ is independent of $\ep$ and is $\tau(\xi)$-periodic in $s$
with zero average, $\int_0^{\tau(\xi)} \gamma (s,\xi) ds = 0$.
%
\begin{equation*}
\begin{aligned}
   \phi (H;\xi_0,\phi_0) &=
   \phi_0 + \int_0^H \left[ \ei g_1 ( \xi (t;\xi_0,\phi_0) ) +g_{2}( \xi (t;\xi_0,\phi_0)) \right] dt
   \\
    &=  \phi_0 + \int_0^H \left[ \ei g_1 (\bar\xi (t) ) +g_{2}(\bar\xi (t)) \right] dt + 
   \int_0^H g'_{1}(\bar\xi (t)) \gamma (t/\ep,\bar\xi(t)) dt + {\mathcal O}(\ep) \\
   &= \phi_0 + F (\xi_0 ; \ep) + {\mathcal O}(\ep) ,
\end{aligned}
\end{equation*}
for some function $F$ that depends only on $\xi_0$ and $\ep$, but not 
on the initial phase $\phi_0$.
In particular, we note that
\begin{equation*}
   \phi (w_1) - \psi_0 = \phi ({\mathcal F}_H w_0 ) - \psi_0 = \phi (H;\xi_0,\psi_0) - \psi_0
      = F(\xi_0;\ep) + {\mathcal O} (\ep) .
\end{equation*}
Similarly,
\begin{equation*}
   \phi_1 - \phi_0 = \phi ({\mathcal F}_H u_0 ) - \phi_0 = \phi (H;\xi_0,\phi_0) - \phi_0
      = F(\xi_0;\ep) + {\mathcal O} (\ep) .
\end{equation*}
Hence, $\phi (w_1) - \psi_0 = \phi_1 - \phi_0 = {\mathcal O}(\ep)$.
In other words, starting at $w_0$ instead of $u_0$ introduces a phase shift that 
is practically constant.
We have then
\begin{equation*}
\begin{aligned}
   t_1^+ &= (\psi_1-\phi_1) \ep\tau_1 
   = \left( \psi_0 - \phi_0 \right) \ep \left[ \tau_0 + \frac{\partial \tau}{\partial \xi}\cdot(\xi_1 - \xi_0) \right]  + {\mathcal O}  (\ep H^2+\ep^2) \\
   &= (\psi_0-\phi_0) \ep \tau_0  + (\psi_0-\phi_0) \ep \frac{\partial \tau}{\partial \xi} \cdot (\xi_1 - \xi_0) + {\mathcal O}  (\ep H^2+\ep^2)\\
&= t_0^+ + H  t_0^+ \Delta + {\mathcal O} (\ep H^2+\ep^2),
\end{aligned}
\end{equation*}
where 
%
\(   
\Delta = \frac{1}{H\tau_0} \frac{\partial \tau}{\partial \xi} \cdot ( \xi_1 - \xi_0) = {\mathcal O} (1).
\)
%
Similarly,
\begin{equation*}
   t_1^- =  t_0^- + H t_0^- \Delta  + {\mathcal O} (\ep H^2+\ep^2).
\end{equation*}
Consider
\begin{equation*}
   {\mathcal F}_{t^\pm_0} u_1 = {\mathcal F}_{- H t^\pm_0 \Delta  } {\mathcal F}_{t_1^\pm} u_1 + {\mathcal O} (\ep H^2+\ep^2) =  {\mathcal F}_{- H t^\pm_0 \Delta } 
   \tilde w_1^\pm + {\mathcal O} (\ep H^2+\ep^2).
\end{equation*}
Expanding around $\tilde w_1^\pm$ 
\begin{equation*}
   {\mathcal F}_{t^\pm_0} u_1 =  
   \tilde w_1^\pm - \frac{H}{\ep} \delta t_0^\pm \Delta + {\mathcal O} (H^2 + \ep^2),
\end{equation*}
for some vector $\delta \in \R^d$ independent of $\ep$.
Therefore, taking a linear combination $\lambda_+ + \lambda_- =1$
and denoting $\lambda=(\lambda_+,\lambda_-)$,
\begin{eqnarray}
   {\mathcal S}_H^\lambda (u_1;u_0,v_0) &=&
   \lambda_+ {\mathcal F}_{t_0^+} u_1 + \lambda_- {\mathcal F}_{t_0^-} u_1\label{eq:forward_alignment_state}  \\
   &=&\left( \lambda_+ \tilde w_1^+ + \lambda_- \tilde w_1^- \right)  
       + \frac{H}{\ep} \delta \left( \lambda_+ t_0^+ + \lambda_- t_0^- \right) \Delta 
  +  {\mathcal O} (H^2+\ep^2) \nonumber\\
   &=& w_1 + \frac{H}{\ep} \delta \left( \lambda_+ t_0^+ + \lambda_- t_0^-  \right) \Delta  
      +  {\mathcal O} (H^2+\ep). \nonumber
\end{eqnarray}
Finally we see that with the choice  
\begin{equation*}
   \lambda_+ = \frac{-t_0^-}{t_0^+ - t_0^-} ~~,~ \lambda_- = \frac{t_0^+}{t_0^+ - t_0^-} ,
\end{equation*}
the first order term cancels.
Thus, we obtain a second order accurate forward alignment to $w_1$. 
See Figure~\ref{fig:Stilde} for the error of \eqref{eq:forward_alignment_state}
for the simple example of the expanding spiral \eqref{eq:simpleSpiral}.



%

Algorithm \ref{alg:pararealFull} applies the convex combination \eqref{eq:forward_alignment_state} 
in forward alignment. Indeed, convergence in the state variable heavily relies on this step because the new point 
$\tilde{u}_{F,n}^{k-1}$ after the forward alignment 
is assigned as the reference for local alignment in the next coarse interval.
See Step 2(c)iii.
We emphasize that taking ${\mathcal S}^{\lambda}_H$ as \eqref{eq:forward_alignment_state} 
may shift the slow coordinates of the resulting $u_1$ from what was computed by the multiscale coarse integrator and assumed to be accurate.
In the next subsection, we  propose a more elaborate method to further improve the overall accuracy of the forward alignment step.

\subsection{Improving accuracy in forward alignment}
Here, the idea is 
that we identify the convex combination with the point  
which divides the trajectory of \eqref{eq:HiOscODEgeneral} 
originating from ${\mathcal F}_{t_0^+} u_1$ and ending close to ${\mathcal F}_{t_0^-} u_1$ by a proportion of $\lambda_{-}$
to $\lambda_{+}$.
Since there are two orientations of ${\mathcal F}_{t} ({\mathcal F}_{t_0^+} u_1)$
defined by forward and backward in time integrations,
the modified convex combination will provide us with two points 
depending on the orientations, and we will choose the one closer to $\mathcal{S}_{H}^{\lambda}u_{1}$.

First, we propose to find the first two local minimizers of 
\begin{equation}
J(t;{\mathcal F}_{t_0^+} u_1, {\mathcal F}_{t_0^-} u_1)=|\mathcal{F}_{t}({\mathcal F}_{t_0^+} u_1)-{\mathcal F}_{t_0^-} u_1|^{2},
\nonumber
\end{equation}
such that
%
$-\ep \tau_1 < \Gamma^-_p < 0 <\Gamma^+_p < \ep \tau_1.$
%
Denoting
\begin{equation}
   t_0^{++} =  t_0^+ + \lambda_- \Gamma^+_p \textrm{~~and~~}t_0^{+-} =   t_0^+ + \lambda_- \Gamma^-_p,
\label{eq:time_pp_pm}
\end{equation}
we again find the first local minimizers of
\begin{equation}
\begin{aligned}
J(t;{\mathcal F}_{t_0^-} u_1, {\mathcal F}_{t_0^{++}} u_1) &= |\mathcal{F}_{t}({\mathcal F}_{t_0^-} u_1)-{\mathcal F}_{t_0^{++}} u_1|^{2},\\
J(t;{\mathcal F}_{t_0^-} u_1, {\mathcal F}_{t_0^{+-}} u_1) &= |\mathcal{F}_{t}({\mathcal F}_{t_0^-} u_1)-{\mathcal F}_{t_0^{+-}} u_1|^{2},
\label{eq:phase_diff0}
\end{aligned}
\end{equation}
such that
%
$-\ep \tau_1 < \Gamma^-_* < 0 \textrm{~~and~~} 0 <\Gamma^+_* < \ep \tau_1,$
%
and denote them by
\begin{equation}
   t_0^{--} =  t_0^- + \Gamma^-_*,\;t_0^{-+} =   t_0^- + \Gamma^+_*.
\label{eq:time_mm_mp}
\end{equation}
With local minimizers of \eqref{eq:phase_diff0},
the phases between 
${\mathcal F}_{t_0^{++}}$ and ${\mathcal F}_{t_0^{--}}$, and between ${\mathcal F}_{t_0^{+-}}$ and ${\mathcal F}_{t_0^{-+}}$
are the same. Now, we define the new weights using $t_0$\rq{}s in \eqref{eq:time_pp_pm} and \eqref{eq:time_mm_mp} by
\begin{equation}
   {\lambda}_{++} = \frac{-t_0^{--}}{t_0^{++} - t_0^{--}} , ~  {\lambda}_{+-} = \frac{-t_0^{-+}}{t_0^{+-} - t_0^{-+}},  ~
   {\lambda}_{--} = \frac{t_0^{++}}{t_0^{++} - t_0^{--}} , ~  {\lambda}_{-+} = \frac{t_0^{+-}}{t_0^{+-} - t_0^{-+}}.
\nonumber
\end{equation}
The convex combination \eqref{eq:forward_alignment_state} is now modified as
\begin{equation}
 {\mathcal S}_H^{{\lambda}_1} (u_1;u_0,v_0) = {\lambda}_{++} {\mathcal F}_{t_0^{++}} u_1 + {\lambda}_{--} {\mathcal F}_{t_0^{--}} u_1,\;\;
 {\mathcal S}_H^{{\lambda}_2} (u_1;u_0,v_0) = {\lambda}_{+-} {\mathcal F}_{t_0^{+-}} u_1 + {\lambda}_{-+} {\mathcal F}_{t_0^{-+}} u_1.
\nonumber
\end{equation}
Here, we note that 
\begin{equation}
 \xi({\mathcal S}_H^{{\lambda}_1} (u_1;u_0,v_0)) = \xi(u_1) + \mathcal{O} (\ep),\;\;
 \xi({\mathcal S}_H^{{\lambda}_2} (u_1;u_0,v_0)) = \xi(u_1) + \mathcal{O} (\ep).
\nonumber
\end{equation}
In words, the modified convex combinations ${\mathcal S}_H^{{\lambda}_1} (u_1;u_0,v_0)$ and
${\mathcal S}_H^{{\lambda}_2} (u_1;u_0,v_0)$ 
guarantee the accuracy of order $\ep$ in the
slow variables of $u_1$.

Now, we propose to implement the forward alignment ${\mathcal{S}}_{H}^{\lambda} (u_1;u_0,v_0)$ as follows.
\medskip

%
%
%

\noindent 
\centerline{
\framebox{\begin{minipage}[c]{1\columnwidth}%
{\it \textbf{Forward alignment: }}
\begin{enumerate}
\item Set the reference point using \eqref{eq:forward_alignment_state},
$\hat{u}:=\lambda_{+}\mathcal{F}_{t_{0}^{+}}u_{1}+\lambda_{-}\mathcal{F}_{t_{0}^{-}}u_{1}$.
\item Compute two modified convex combinations with opposite orientations,
\begin{equation}
{\lambda}_{++} {{\mathcal F}}_{t_0^{++}} u_1 + {\lambda}_{--} {{\mathcal F}}_{t_0^{--}} u_1,\;
{\lambda}_{+-} {{\mathcal F}}_{t_0^{+-}} u_1 + {\lambda}_{-+} {{\mathcal F}}_{t_0^{-+}} u_1.
\nonumber
\end{equation}
\item Denote by $\mathcal{S}_{H}^{\lambda} (u_1;u_0,v_0)$ the combination closer to $\hat{u}$.
\end{enumerate}
\end{minipage}}
}
\medskip
\begin{rem}
An unperturbed system of the HiOsc system \eqref{eq:HiOscODEgeneral}, if exists, preserves the slow variables
but changes the fast variables. Indeed,
by denoting $\mathcal{F}^{0}$ the fine integrator for the unperturbed system, $\xi(\mathcal{F}_t^{0} u_0)=\xi(u_0)$ for all $t>0$ but $\phi(\mathcal{F}_t^{0} u_0) \neq \phi(u_0)$ for some $t>0$. 
If the unperturbed system of \eqref{eq:HiOscODEgeneral} 
is explicitly known, one can achieve more accurate local and forward alignments by using $\mathcal{F}^{0}$ in the minimization of $J(t)$.
Unfortunately, the unperturbed system is not explicitly known for general systems.
Comparison of the parareal solutions using $\mathcal{F}$  with $\mathcal{F}^{0}$ 
will be presented in Sections~\ref{sec:spiral1} and \ref{sub:FPU}.
\end{rem}

\section{A multiscale integrator based on Poincar\'e-map}
\label{sec:reviewPoincare}
\setcounter{equation}{0}

Even though the goal of the multiscale system is a consistent description
of only the slow variables, in practice,  obtaining an explicit expression 
for the slow variables is often difficult or impossible,
in particular for high-dimensional systems
(see \cite{Artstein-Kevrekidis-Slemrod-Titi11} for an example).
In \cite{BFHMM2012}, a new type of multiscale methods
using a Poincar\'e map technique was introduced. This method 
only assumes the existence of slow variables but does not 
use its explicit form.
A novel on-the-fly filtering technique achieves high order accuracy.
Recall the general two-scale ODE \eqref{eq:HiOscODEgeneral}
with initial condition $u_0 \in D \subset\mathbb{R}^{d}$:
\begin{equation}
   \dot{u}=\epsilon^{-1}f_{1}(u)+f_{0}(u),\,\,\, u(0)=u_0 .
   \label{eq:singularODE3}
\end{equation}
By ignoring the lower order perturbation part of the vector field, an unperturbed dynamical system is defined.
The essential part of the Poincar\'e-map technique is to generate a path whose projection on the slow subspace has the correct slow dynamics.
To this end, the scheme solves both the perturbed and the unperturbed systems 
from the same initial conditions for short time intervals, 
and compares the resulting trajectories.

The method relies on the following assumptions regarding the HiOsc dynamics
\begin{assumption}
\label{assumption1}
The dynamics of the unperturbed equation
\begin{equation}
\dot{v} =\ei f _1 (v),\,\,\, v(0)=v_0.
\label{eq:msFastOnly}
\end{equation}
is ergodic with respect to an invariant manifold $\mathcal{M}(v_{0})$.
\end{assumption}
\noindent
We denote the solution of \eqref{eq:msFastOnly} by $v(t;v_0)$. 
\begin{assumption}
\label{assumption2}
The invariant manifolds $\mathcal{M}(z)$
is defined by the intersection of the level sets of slow
variables $\xi_{1},\xi_{2},\cdots,\xi_{k},$ $k<d$. 
More precisely, we may identify the invariant manifold of $v$
by level sets of the slow variables for u, 
$\mathcal{M}(z)=\cap_{j=1}^{k}\{\zeta\in\mathbb{R}^{d}:\xi_{j}(z)=\xi_{j} (\zeta) \}$.
\end{assumption}
\noindent
Hence, the solution $u(t)$ defines a foliation of invariant manifolds
$\mathcal{M}(t):=\mathcal{M}(u(t))$.
Note that our method only assumes the existence of such $\xi$'s
but does not require obtaining them.

Suppose we solve the full equation \eqref{eq:singularODE3} and 
the associated unperturbed version \eqref{eq:msFastOnly}
with the same initial condition.
Then, it is possible to extract
the flow of $\mathcal{M}(t)$ from comparison of $u(t)$
and $v(t)$ without explicitly knowing the slow variables. 
The central idea is to locally create a path $\gamma$ in states space 
that is transversal to the fast flow.
This cut will be defined and approximated by a procedure that 
realizes a Poincar\'e return map along it.
We shall look for a slow $\gamma(t)$, i.e., require that $| \dot{\gamma} |\le C$
such that for any slow variable $\xi$, $\xi (\gamma(t))=\xi (u(t))$. 
In other words, the effective slow path $\gamma (t)$ goes through the same foliation
of slow manifolds as the exact solution,
$\mathcal{M}(\gamma(t))) = \mathcal{M}(u(t))$.
The time derivatives of such effective paths can be obtained by extracting
the influence of lower order perturbations in the given oscillatory equation.
Approximating the derivative will require solving the HiOsc system
for reduced time segments of order $\ep$.
Since $\gamma$ is slow, it can be approximated using macroscopic
integrators with step size $H$ independent of $\ep$.
As a result, the overall computational complexity of the resulting 
algorithm is sublinear in $\epsilon^{-1}$.

\begin{figure}[H]
   \begin{centering}
   \includegraphics[width=10cm]{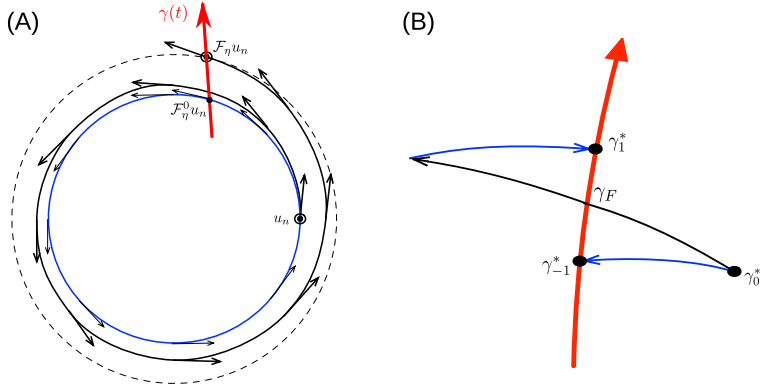} 
   \par\end{centering}
   \caption{\label{fig:BFHMM_schematics} The Poincar\'e map-type technique approximated all slow variables but does not require knowing their explicit formulas. (A) A forward-Euler type
   construction. (B) The symmetric Poincar\'e method is first order, but symmetric
   with respect to the Poincar\'e return points.}
\end{figure}

To be consistent with previous notation, we
denote by ${\mathcal F}_t$ the fine scale approximated propagator
for the full equation \eqref{eq:HiOscODEgeneral}, and 
by ${\mathcal F}_t^0$ the fine scale approximated propagator
for the unperturbed equation \eqref{eq:msFastOnly}
in which the low order perturbation is turned off.
In particular, note that under
the dynamics of \eqref{eq:msFastOnly} all slow variables are constants of motion.
Let $\ep<\eta<H$. 
A basic Forward-Euler step, depicted in Figure~\ref{fig:BFHMM_schematics}A  can be written as
\begin{equation}
   u_{n+1} = {\mathcal F}_\eta u_n + \frac{H}{\eta} \left( {\mathcal F}_\eta u_n - {\mathcal F}^0_\eta u_n  \right) .
\label{basic_FE}
\end{equation}
The values of the effective path $\gamma(t)$ at $t=nH$ is then identified with $u_n$, $n=0,1,\cdots,N$.
High order approximations of $\gamma (t)$ may be obtained by
combining several steps and using high-order extrapolation.
The name "Poincar\'e technique" alludes to the fact that $\gamma (t)$ is
transversal to the solution curves of the unperturbed equation. 
Thus, the full solution induces a Poincar\'e return map,
which is used to approximate $\gamma (t)$.
See \cite{BFHMM2012} for details.

The bottleneck in the efficiency of
the new algorithm is a consequence of small-amplitude, high-frequency oscillations
in $\xi (u(t))$.
The accuracy can be improved by sampling the derivatives of a locally smooth
average of $\xi$, $\bar{\xi}$ instead of the weakly oscillating $\xi$. 
Since we assume no explicit knowledge
about the slow variables, $\bar{\xi}$ must be computed intrinsically.
In \cite{BFHMM2012}, a filtering technique is
proposed for the simple case in which the invariant manifold of
the unperturbed equation is diffeomorphic to a circle, i.e., the unperturbed 
dynamics is periodic.
More precisely, we propose to replace \eqref{eq:HiOscODEgeneral}
by the filtered equation
\begin{equation}
   \dot{\bar{u}}=\ei f_1(\bar{u})+K (t/\eta) f_0 (\bar{u},t,\ei t),\,\,\,\, t^{*}\le t\le t^{*}+\eta,
   \label{eq:filtered-eqn}
\end{equation}
where the filter $K_\eta$ is $C^q ([0,\eta])$ and satisfies the moment condition of the form
\begin{equation}
   \int_{0}^{\eta}K_{\eta}(\eta-s)s^{j}ds=\int_{0}^{\eta}s^{j}ds,j=0,1,2,\cdots,p.
 \nonumber
\end{equation}
%
%
%
The accuracy of \eqref{eq:filtered-eqn} was demonstrated and analyzed in \cite{BFHMM2012}.


%

\subsection{Symmetric Poincar\'e methods}
\label{sec:z}

The simple Forward-Euler step \eqref{basic_FE} can be applied in
simple situations in which the frequency of the fast oscillation is not a slow variable itself 
(i.e., $g_1$ in \eqref{eq:odefastslow} is not a function of $\xi$) \cite{BFHMM2012}.
This restriction can be lifted by generating interpolation points $\gamma^{*}_k$ symmetrically
described as follows.
Our idea it to generate and choose interpolation points $\gamma^{*}$ 
so that in the state space, 
$(2\eta)^{-1} \left( \gamma_1^* - \gamma_{-1}^* \right)$
approximates \textit{implicitly} the derivative of $\bar{\xi}(t)$
but results in small derivative of $\phi(t)$,
more precisely, of order $\eta^2$.
%
%
The method, originally proposed and analyzed in \cite{SeongJunThesis},
can be described as
\begin{equation}
   u_{n+1} = \gamma_{-1}^* + \frac{H}{2\eta} \left( \gamma_1^* - \gamma_{-1}^* \right) ,
\label{eq:PoincareRule}
\end{equation}
where
\begin{equation}
   \gamma_{-1}^* = {\mathcal F}^0_{\eta} u_n, \;\;
   \gamma_{1}^* = {\mathcal F}^0_{-\eta} {\mathcal F}_{ 2\eta} u_n .
\nonumber
\end{equation}
%
Convergence of the method is proved in Appendix B.
See also \cite{SeongJunThesis} for different types of Poincar\'e methods.
The method \eqref{eq:PoincareRule} defines a propagator, denoted ${\mathcal M}_H$,
\begin{equation}
   {\mathcal M}_H =  \left( 1 - \frac{H}{2\eta} \right) {\mathcal F}^0_{\eta} + \frac{H}{2\eta}
   {\mathcal F}^0_{-\eta} {\mathcal F}_{2 \eta} .
\nonumber
\end{equation}
The efficiency of the combined parareal with multiscale Poincar\'e method can be evaluated by counting the number of fine steps in
each parareal iteration.
An effective time segment $\tau$ is denoted the computational cost of a single parareal 
iteration equivalent to fine integration of length $\tau$:
\begin{equation}
   \tau= H + h_{fine} \eta \left( \frac{4}{h_{Poincare}} + \frac{1}{h_{phase}} \right) \left( 1 + \frac{T}{H} \right)
\label{eq:cost1}
\end{equation}
%
where $h_{fine}$, $h_{Poincare}$, and $h_{phase}$ are
the microscopic step sizes used for
fine, Poincar\'e and phase alignment methods, respectively.
Overall, with $K$ iterations, the computational speed-up (assuming maximal parallelization)
compared to 
direct numerical simulation is $T/K \tau$.

\section{Numerical examples}
\label{sec:examples}
\setcounter{equation}{0}

\subsection{Expanding spiral I}
\label{sec:spiral1}

Consider the following nonlinear equation in the complex plane
\begin{equation}
\dot{u}=i\epsilon^{-1}u|u|+u|u|^{-1},\label{eq:num_spiral-1}
\end{equation}
with the initial value $u(0)=1$. 
The associated unperturbed equation is known to be
\begin{equation}
\dot{v}=i\epsilon^{-1}v|v|,\;\; v(0)=v_0. \label{eq:num_spiral1-unper}
\end{equation}
As in (\ref{eq:odefastslow}),
the dynamics of $u(t)$ can be analyzed by the corresponding system
of slow and fast variables:
\begin{equation}
	\begin{cases}
	\dot{\xi}=1, & \xi(0)=1,\\
	\dot{\phi}=\epsilon^{-1}\xi, & \phi(0)=0.
	\end{cases}\label{eq:num_spiral_diffeo-1}
\end{equation}
We see immediately from Definition\textbf{ \ref{def:slowvariable}}
that $\xi$ is a slow variable. The difficulty in the phase alignment
lies in the singular term in the RHS of $\dot{\phi}$. 
Note that (\ref{eq:num_spiral_diffeo-1})
is never used in our algorithm as $\xi$ in (\ref{eq:num_spiral_diffeo-1})
is only used to show the convergence in the slow variables.

\begin{figure}[H]
\begin{centering}
\subfloat[$\epsilon=10^{-3}$]{\begin{centering}
\includegraphics[width=5.5cm]{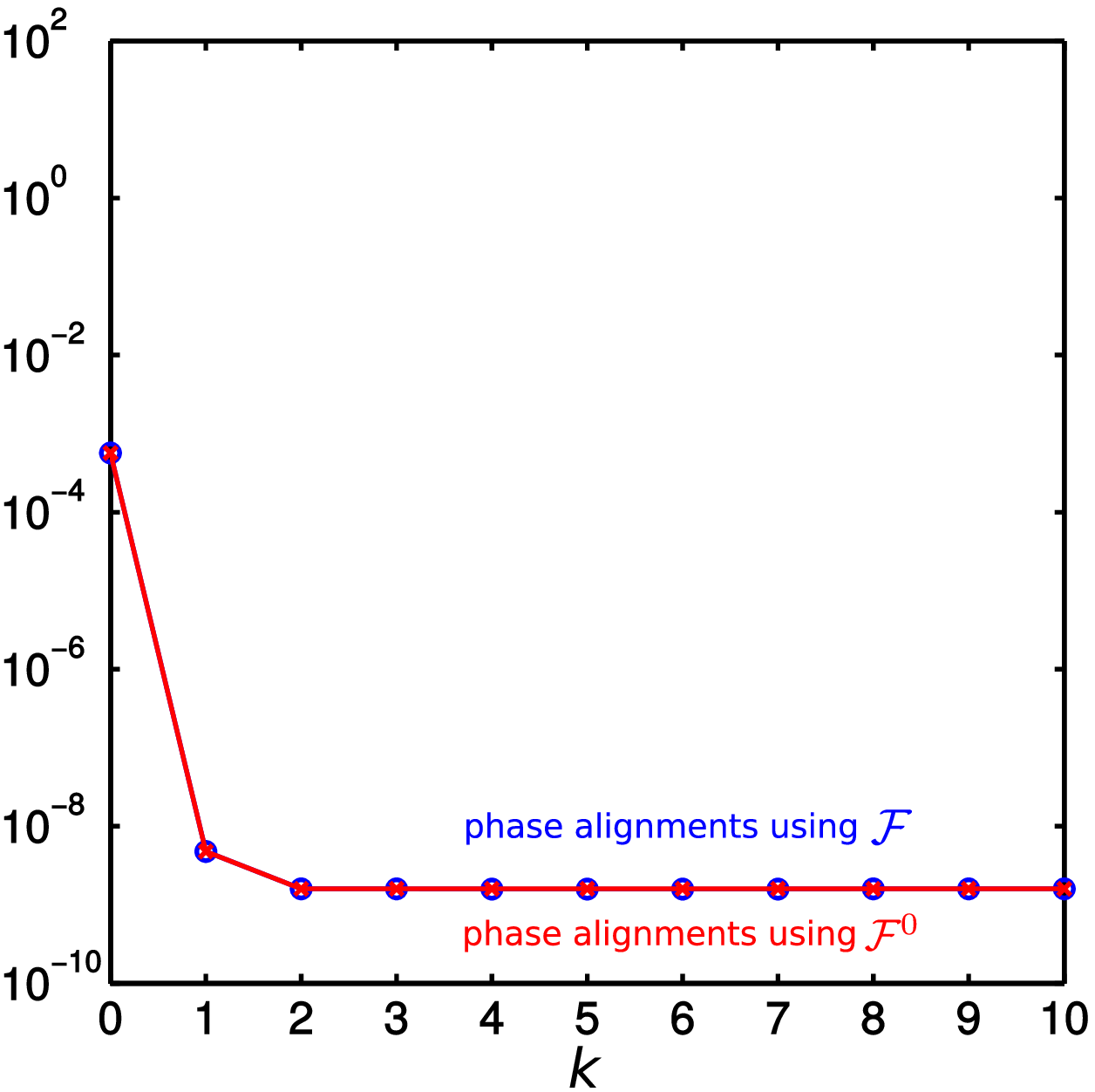}
\includegraphics[width=5.5cm]{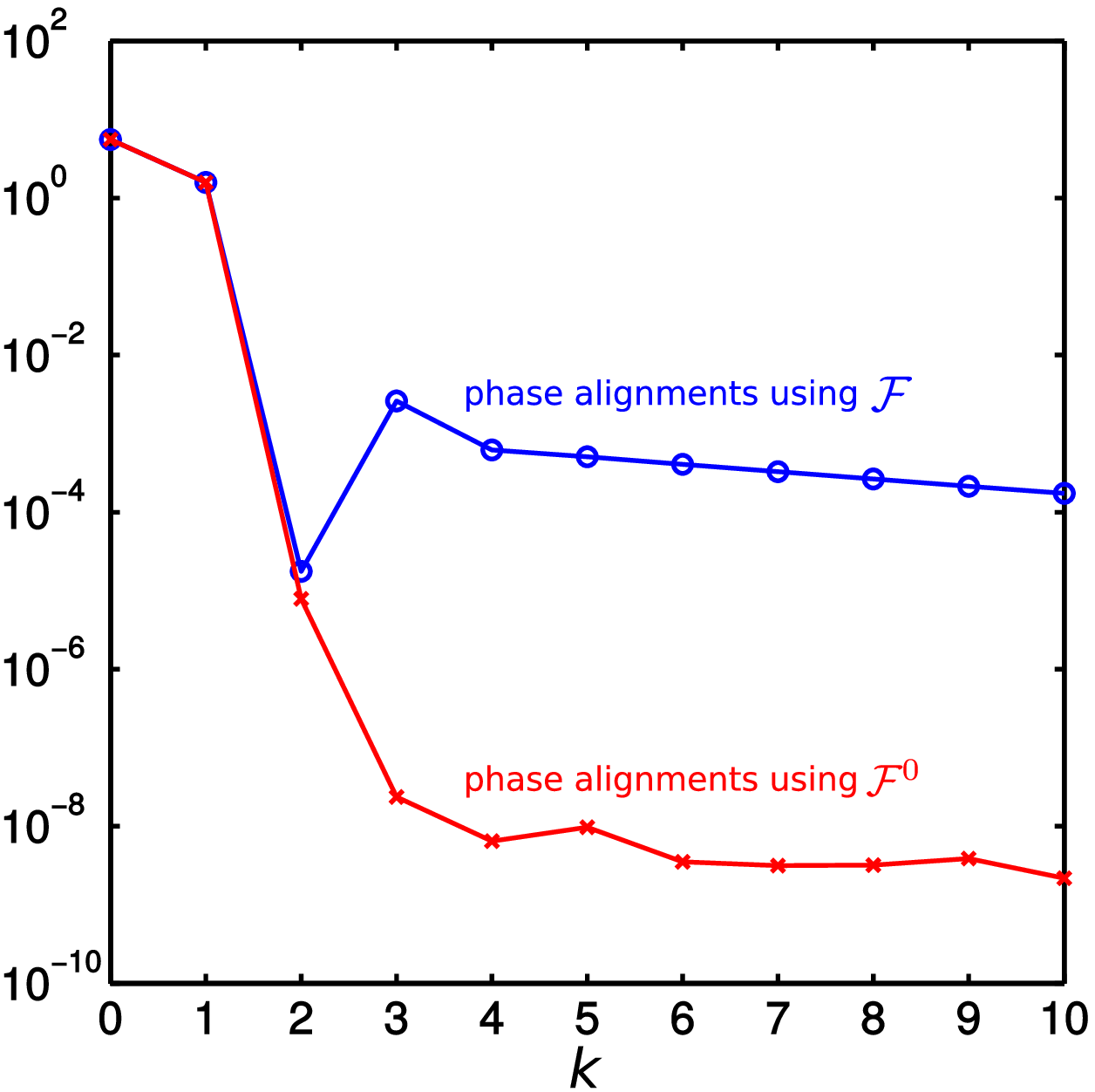}
\par\end{centering}
}
\par\end{centering}
\vspace*{-0.2cm}
\begin{centering}
\subfloat[$\epsilon=10^{-4}$]{\begin{centering}
\includegraphics[width=5.5cm]{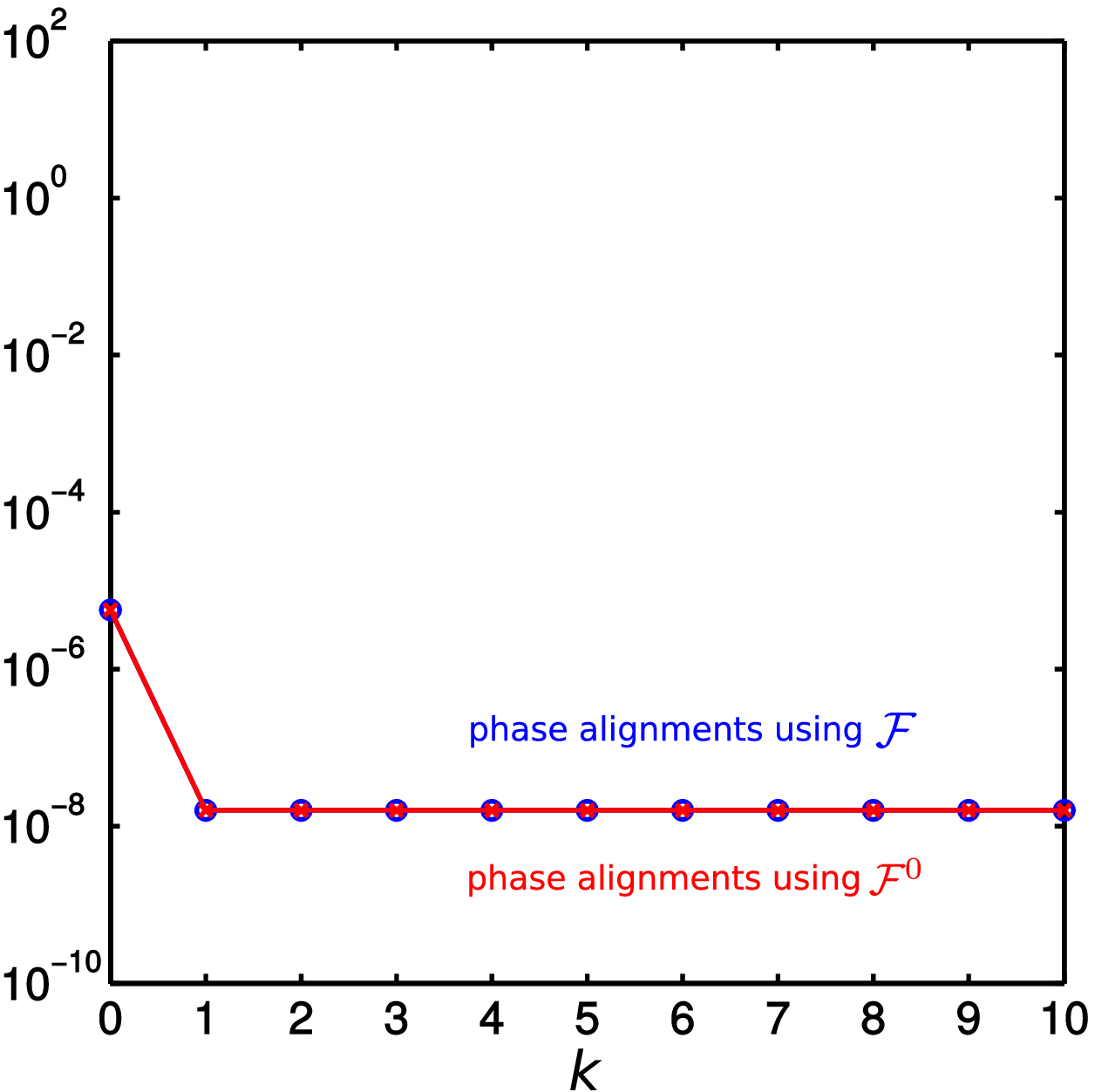}
\includegraphics[width=5.5cm]{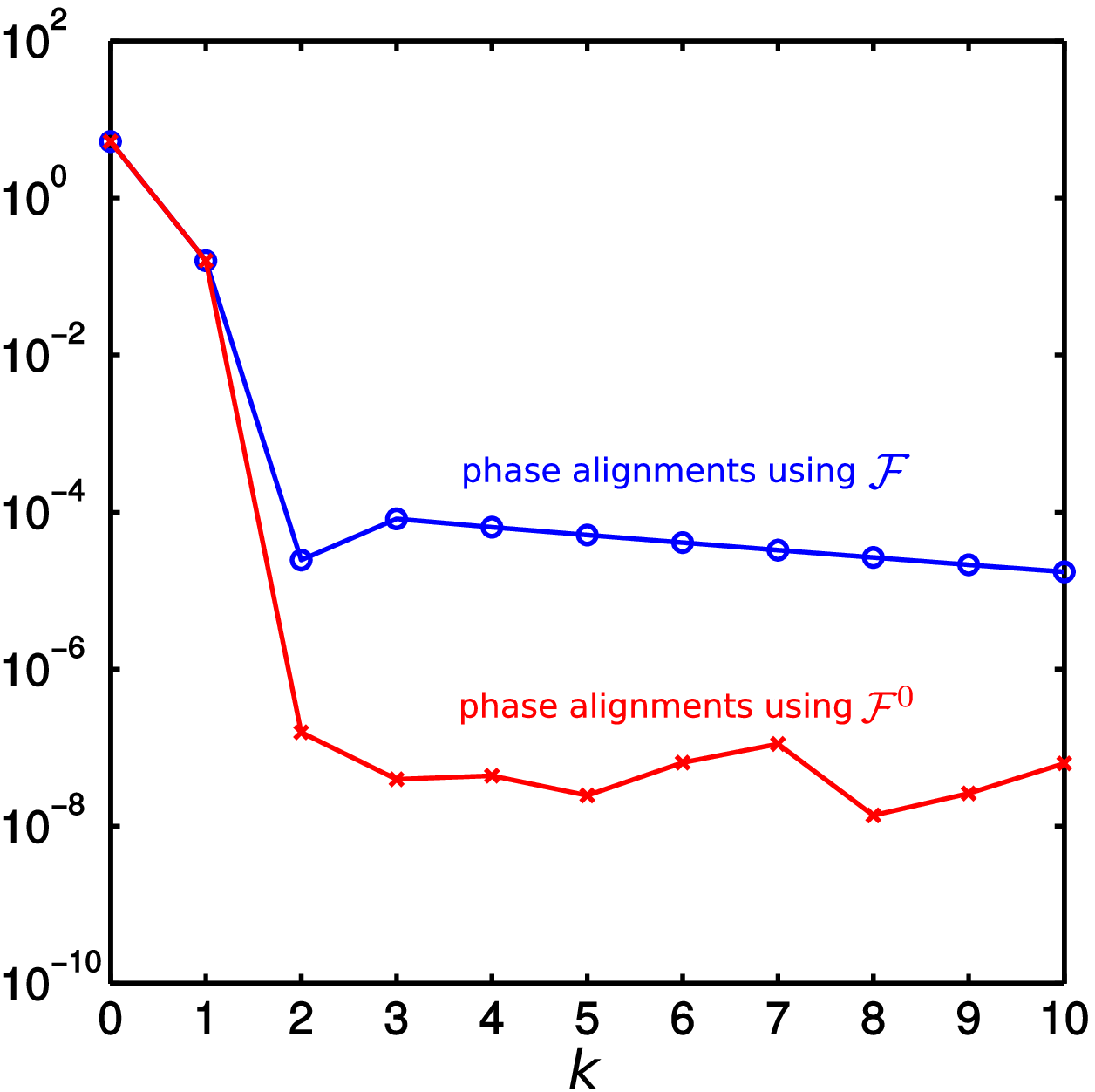}
\par\end{centering}
}
\par\end{centering}
\vspace*{-0.2cm}
\caption{Expanding spiral I, Example~\ref{sec:spiral1}. (Left) Plot of 
$||\xi (\cdot)-\xi \circ u^{k}(\cdot)||_{L^{\infty}([0,T])}$ as a function of iteration
and (Right) the absolute error in the state variables as a function of iteration.
\label{fig:spiral results-1}}
\end{figure}

In this example, we use Algorithm \ref{alg:pararealFull}, ODE45 as a fine integrator and 
the Poincar\'e 2nd order multiscale method
as a coarse integrator (the Midpoint rule
macro-solver and ODE45 micro-solver with z-shape construction of $\gamma$)
to compute the solution. In each micro-simulation, we solve the filtered equation
\begin{equation}
\dot{u}_{n}(t)=i\epsilon^{-1}u_{n}(t)|u_{n}(t)|+K_{\eta}(t-t_{n}) u_{n}(t)|u_{n}(t)|^{-1} ,\,\,\, t_{n}\le t\le t_{n}+\eta,
\label{eq:poincare}
\end{equation}
with a $C^4$ kernel with $p=1$ supported on $[0,\eta]$
is used.
 The parameters are specified in Table \ref{tab:parareal_parameter-spiral6_1-1}
where $\eta_{Poincare}$ and $h_{Poincare}$ are parameters for \eqref{eq:poincare}.
The absolute errors in the slow variables and in the state variables 
as a function of parareal iterations for each
different $\epsilon=10^{-3}$ and $10^{-4}$ are shown in Figure
\ref{fig:spiral results-1}.  
In addition, we compare the errors when both local and forward alignments are established using the full system 
\eqref{eq:num_spiral-1} and the unperturbed system \eqref{eq:num_spiral1-unper}.
The parareal solutions using \eqref{eq:num_spiral-1}
in phase alignments have an error of $\mathcal{O} (\ep)$ in the state variables.
On the other hand, the phase alignment using \eqref{eq:num_spiral1-unper} shows an error in the state variables comparable to that in the slow variables. 
Indeed, knowing the unperturbed equation \eqref{eq:num_spiral1-unper} is an advantage.
In applications, obtaining an explicit expression for the unperturbed equation may not be possible, 
particularly for nonlinear systems.
\begin{table}[H]
\centering{}\caption{Parareal parameters in Example~\ref{sec:spiral1}. \label{tab:parareal_parameter-spiral6_1-1}}
\begin{tabular}{|c|c|c|c|c|c|}
\hline 
$T$ & $H$ & $h_{fine}$ & $\eta_{Poincare}$ & $h_{Poincare}$ & RelTol, AbsTol (ODE45 parameters)\tabularnewline
\hline 
\hline 
2 & 1/10 & $\epsilon/200$ & 7$\epsilon$ & $\epsilon/10$ & $10^{-13}$, $10^{-11}$\tabularnewline
\hline 
\end{tabular} 
\end{table}

\subsection{Expanding spiral with slowly varying fast oscillations}
\label{sec:spiral}

Consider the following HiOsc example
\begin{equation}
\begin{aligned}
   & \dot{x} =   - 2 \pi \ei  [1+ (1-a z_1 ) z_2 ] y + b x \\
   & \dot{y} = 2 \pi \ei [1+(1-a z_1 ) z_2 ] x + b y  \\
   & \dot{z}_1 =1  \\
   & \dot{z}_2 = -a z_2 
   \label{eq:spiral}
\end{aligned}
\end{equation}
where $a,b>0$ are constants.
Initial conditions are $(x,y,z_1,z_2)(0)=(1,0,0,1)$.
The solution of \eqref{eq:spiral} is 
$x(t)= e^{bt} \cos \left[ 2 \pi \ei (1+e^{-at}) t \right]$, $y(t)=e^{bt} \sin \left[ 2 \pi \ei (1+e^{-at}) t \right]$,
$z_1 (t)=t$, and $z_2 (t)=e^{-at}$.

Hence $I=x^2+y^2$, $z_1$ and $z_2$ are three slow variables
while $(x,y)$ is a linear oscillator with expanding amplitude $\sqrt{I}$
and a slowly changing period $\epsilon/(1+e^{-at})$.
The example falls under the general category of HiOsc systems
in which the dynamics of the fast phase slowly evolves according to the slow variables.
The main difference between this example and Example~\ref{sec:spiral1}
is that the derivative of slow variables is not a constant.
As a result, the local ${\mathcal O}(H^2)$ error introduced by the 1st order
coarse multiscale integrator is realized.

The system \eqref{eq:spiral} is integrated using 
the full multiscale Poincar\'e-parareal method,
applying the corrected phase shift described in Section~\ref{seq:HiOscParareal}
to ensure convergence in the state variable.
We stress that the numerical approximation is obtained without using our knowledge
that the system can be decomposed into the three slow variables $I$,
$z_1$ and $z_2$ and a fast phase-like variable $\phi=\arctan (y/x)$.
This decomposition and the exact solution are only used
in order to explain the fast-slow structure in the dynamics
and for demonstrating the rate of convergence of different variables.

Figure~\ref{fig:spiralII}A shows the error in the slow variables
as a function of iteration.
After a single iteration the error in the slow variables drops below~$\ep$,
which is the theoretical limit possible with multiscale methods
on their own.
Figure~\ref{fig:spiralII}B shows the absolute error in the state variable of the entire trajectory.
Initially, the absolute error is large.
This is because the inaccurate slow variables create a jump in the phase 
between coarse time segments.
However, after two iterations, the phase shift becomes accurate and 
the method converges to the exact solution.
Parameters are detailed in Table~\ref{tab:parareal-parameter-spiral6_2}.

\begin{table}[H]
\centering{}\caption{Parareal parameters in Example~\ref{sec:spiral}. \label{tab:parareal-parameter-spiral6_2}}
\begin{tabular}{|c|c|c|c|c|c|c|}
\hline 
$\epsilon$ & $T$ & $H$ & $h_{fine}$ & $\eta_{Poincare}$ & $h_{Poincare}$ & RelTol, AbsTol(ODE45 parameters)\tabularnewline
\hline 
\hline 
$10^{-3}$ & 2 & 0.1 & $\ep/200$ & 7$\epsilon$ & $\epsilon/10$ & $10^{-13}$, $10^{-11}$\tabularnewline
\hline 
\end{tabular} 
\end{table}

\begin{figure}[tbh]
\begin{centering}
\includegraphics[width=5cm]{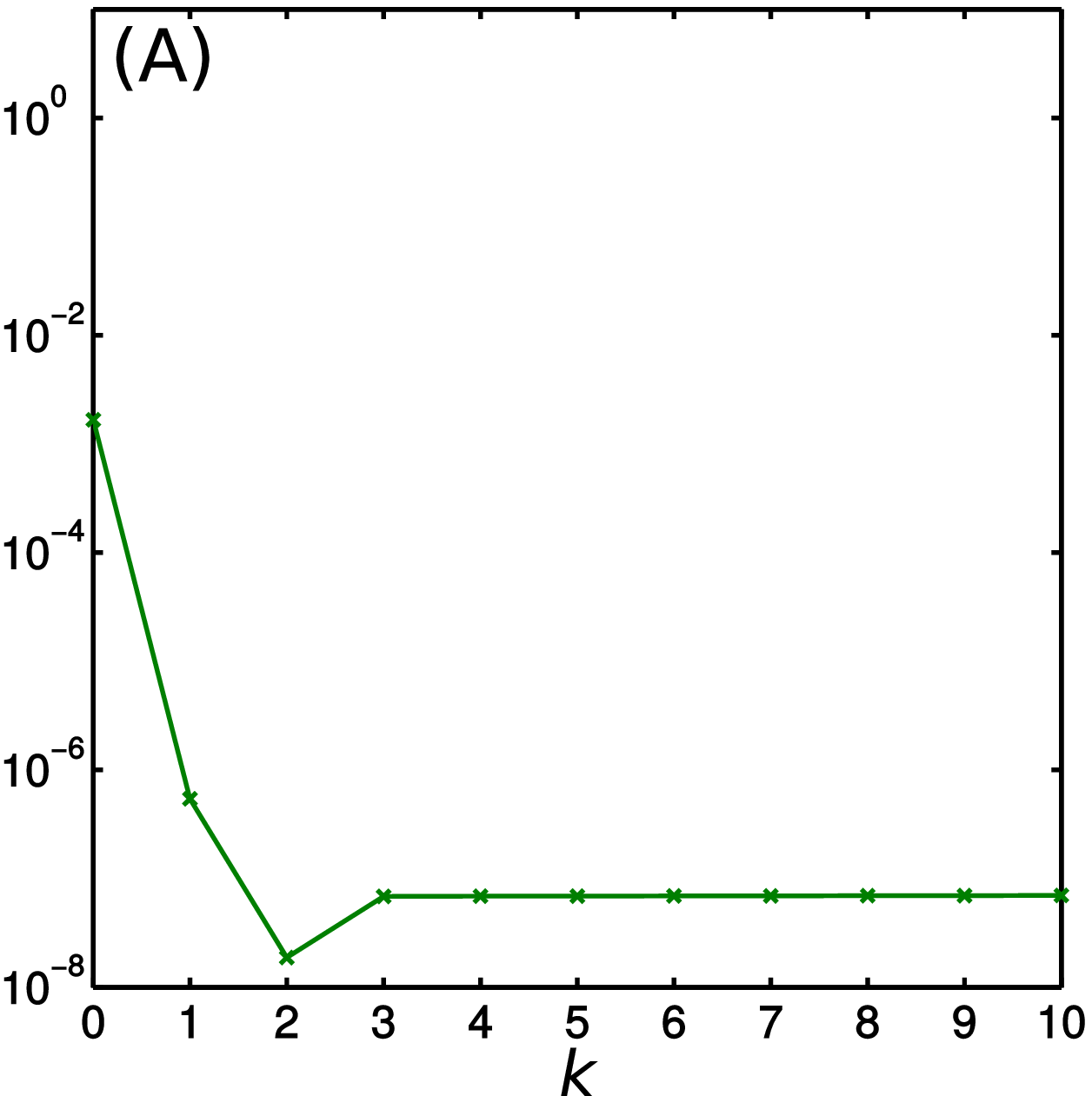} \;
\includegraphics[width=5.01cm]{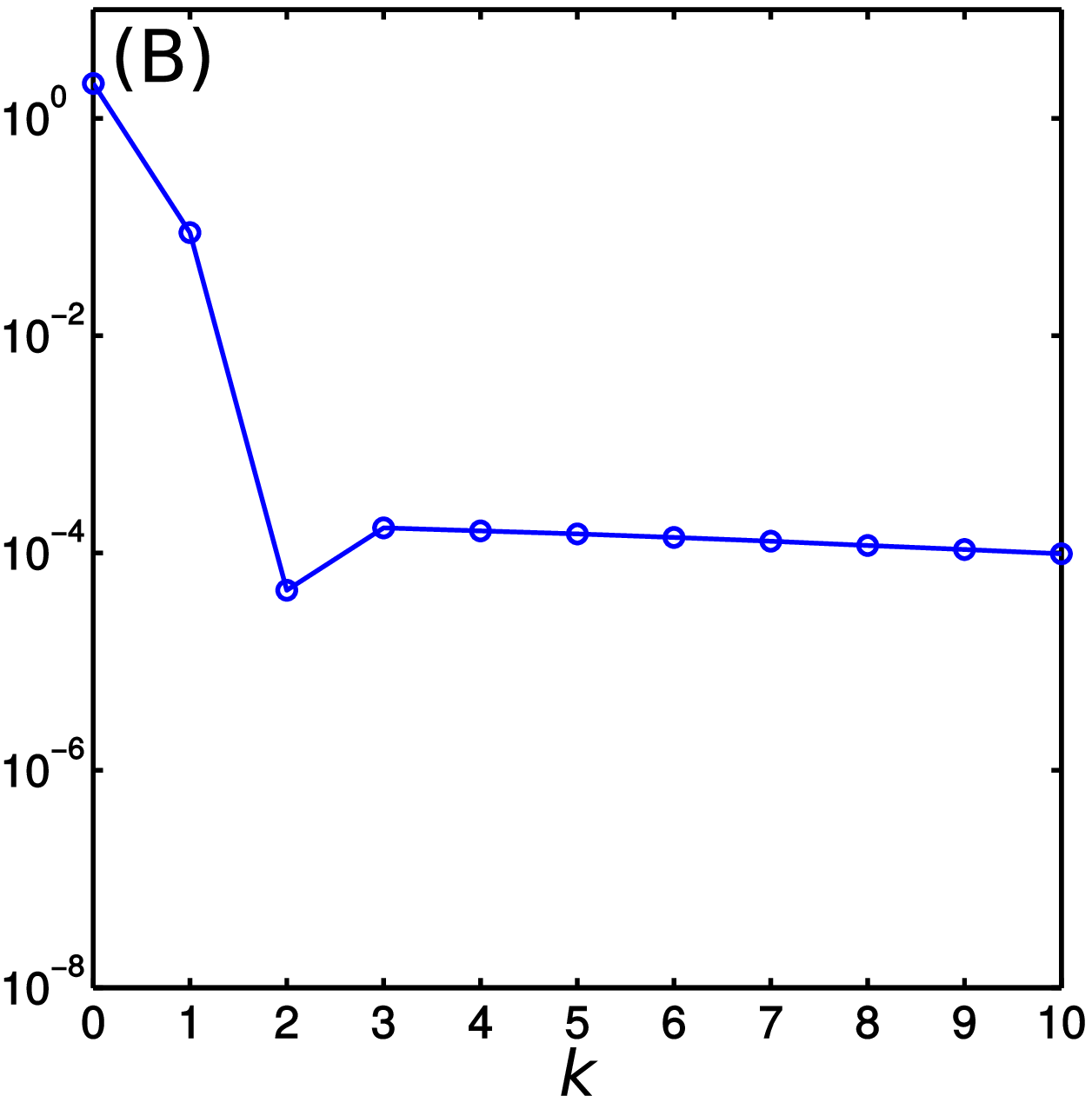}
\par\end{centering}
\caption{Expanding spiral  with $a=0.2$, $b=0.1$, Example~\ref{sec:spiral}. 
(A) The error in the slow variables, $\underset{i=1,2,3}{\max}||\xi_{i}(\cdot)-\xi_{i}\circ u^{k}(\cdot)||_{L^{\infty}([0,T])}$, as a function of iteration.
(B) The absolute error in the state variables as a function of iteration.
\label{fig:spiralII}}
\end{figure}

\subsection{Stellar orbits in a galaxy}
\label{sec:stellar}
The following system is taken from the theory of stellar orbits in a galaxy \cite{Kevorkian-Cole:Perturbations1,Kevorkian-Cole:Perturbations2}:

\[
\begin{cases}
\ddot{r}_{1}+a^{2}r_{1} =\epsilon r_{2}^{2},\\
\ddot{r}_{2}+b^{2}r_{2} =2\epsilon r_{1}r_{2}. 
\end{cases}
\]
Following a change of variables 
$\mathbf{x}=[x_1, v_1, x_2, v_2]^T=[r_1, \frac{d}{d\tau}r_1/a, r_2, \frac{d}{d\tau} r_2/b]^T$
and after a rescaling of time, $t=\epsilon \tau$,
the system can be written in the following form
\begin{equation}
\begin{aligned}
   & \dot{x}_1 =  \ei a v_1 \\
   & \dot{v}_1 = - \ei a x_1 + x^2_2/a \\
   & \dot{x}_2 =  \ei v_2  \\
   & \dot{v}_2 = -\ei x_2 + 2x_1 x_2 /b .
   \label{eq:stellar}
\end{aligned}
\end{equation}
Initial conditions are  $(x_1,v_1,x_2,v_2)(0)=(1,0,1,0)$. Resonance of oscillatory modes take effect in the lower order term when $a=2$ and $b=1$. Using the algorithm proposed in \cite{Ariel-Engquist-Tsai07}, three independent slow variables are identified as
\begin{equation}
\xi_{1}=x_{1}^{2}+v_{1}^{2},\,\,\,\,\,\,\,\,\xi_{2}=x_{2}^{2}+v_{2}^{2},\,\,\,\,\,\,\,\,\xi_{3}=x_{1}x_{2}^{2}+2v_{1}x_{2}v_{2}-x_{1}v_{2}^{2}.\label{eq:stellar-slow-vars}
\end{equation}
The example falls under the category of HiOsc systems
in which two stiff harmonic oscillators are coupled.
The local ${\mathcal O}(H^2)$ error introduced by the 1st order coarse multiscale integrator is realized.

The system \eqref{eq:stellar} is integrated using 
Algorithm~\ref{alg:pararealFull}, 
applying the global alignment algorithm described in Section~\ref{seq:HiOscParareal}
to ensure convergence in the state variable.
The Poincar\'e method is used as a coarse solver, using the trajectory of \eqref{eq:stellar} as the flow ${\mathcal F}$ 
and the one of \eqref{eq:stellar} without a lower order perturbation as the flow ${\mathcal F}^0$.
We stress that the numerical approximation is obtained without using our knowledge
that the system can be decomposed into the three slow variables $\xi_1$,
$\xi_2$ and $\xi_3$ and a fast phase-like variable $\phi$. 
This decomposition are only used
in order to explain the fast-slow structure in the dynamics.

Figure~\ref{fig:spiralII} shows the absolute error in the state variable of the entire trajectory.
Initially, the absolute error is large because the inaccurate slow variables create a jump in the phase 
between coarse time segments.
However, after four iterations the error in the state variable drops below~$\ep$,
which is the theoretical limit possible with multiscale methods
on their own.  
Parameters are detailed in Table~\ref{tab:parareal-parameter-spiral6_2}. The fine integrator is ODE45 method, and 
the coarse integrator is the Poincar\'e 2nd order multiscale method.
A $C^3$ kernel with $p=1$ is used for the filtered equation.

\begin{table}[H]
\centering{}\caption{Parareal parameters in Example~\ref{sec:stellar}. \label{tab:parareal-parameter-spiral6_2}}
\begin{tabular}{|c|c|c|c|c|c|c|}
\hline 
$\epsilon$ & $T$ & $H$ & $h_{fine}$ & $\eta_{Poincare}$ & $h_{Poincare}$ & RelTol, AbsTol(ODE45 parameters)\tabularnewline
\hline 
\hline 
$10^{-4}$ & 14 & 0.5 & $\ep/100$ & 20$\epsilon$ & $\epsilon/10$ & $10^{-13}$, $10^{-11}$\tabularnewline
\hline 
\end{tabular} 
\end{table}

\begin{figure}[tbh]
\begin{centering}
\includegraphics[width=5cm]{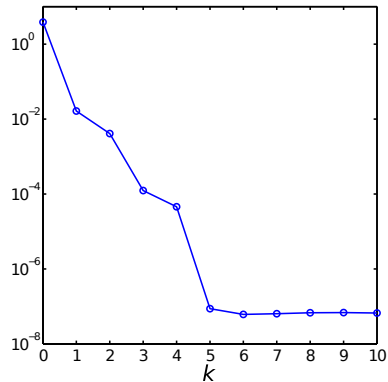}
\par\end{centering}
\caption{Stellar orbits in a galaxy, Example~\ref{sec:stellar}. 
The absolute error in the state variables as a function of iteration is depicted.
\label{fig:spiralII}}
\end{figure}

\subsection{Non-linear oscillators}
\label{sec:VL_oscillator}
\begin{figure}[tbh]
\begin{centering}
\includegraphics[width=5.05cm]{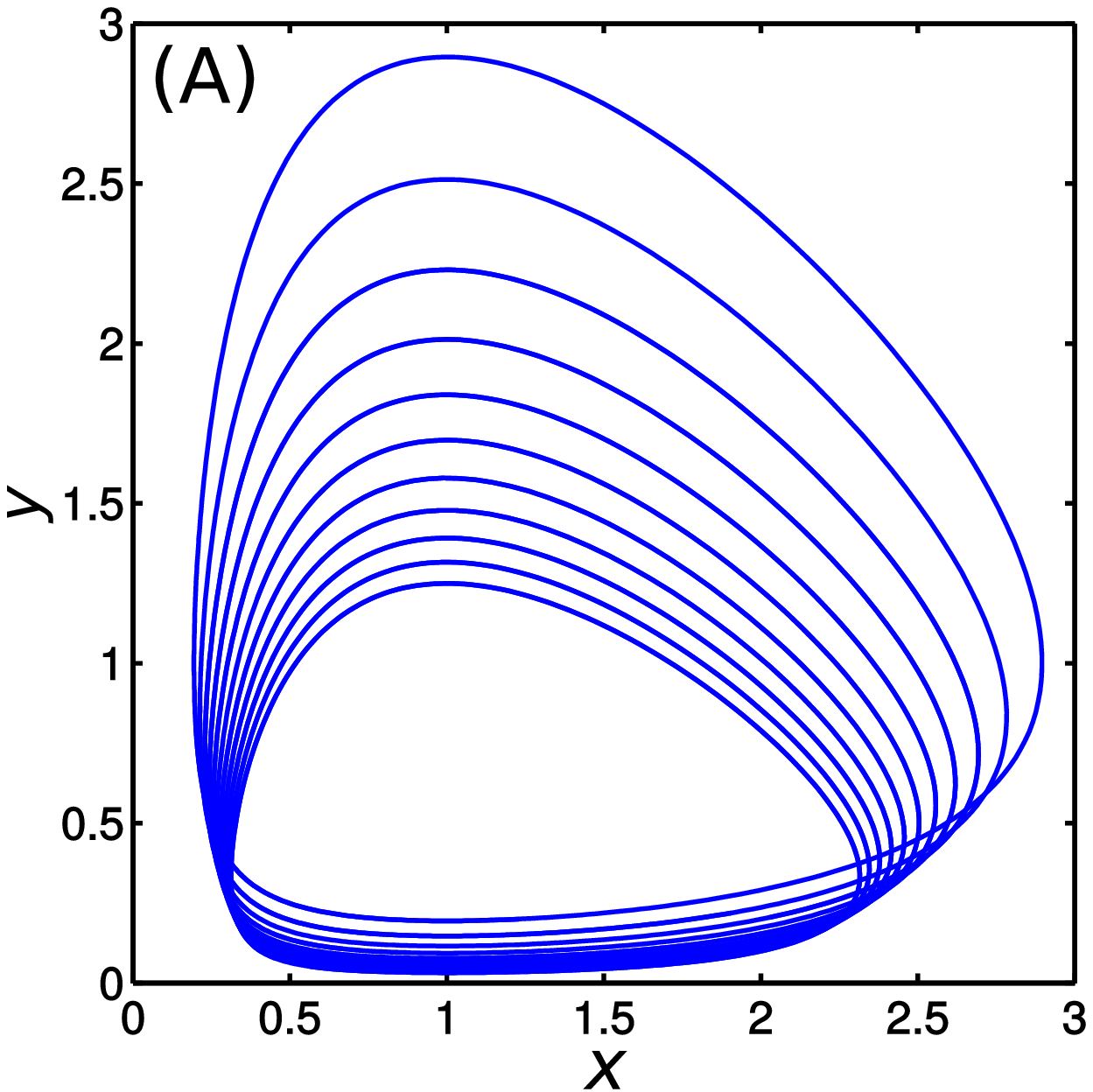} \;
\includegraphics[width=5.1cm]{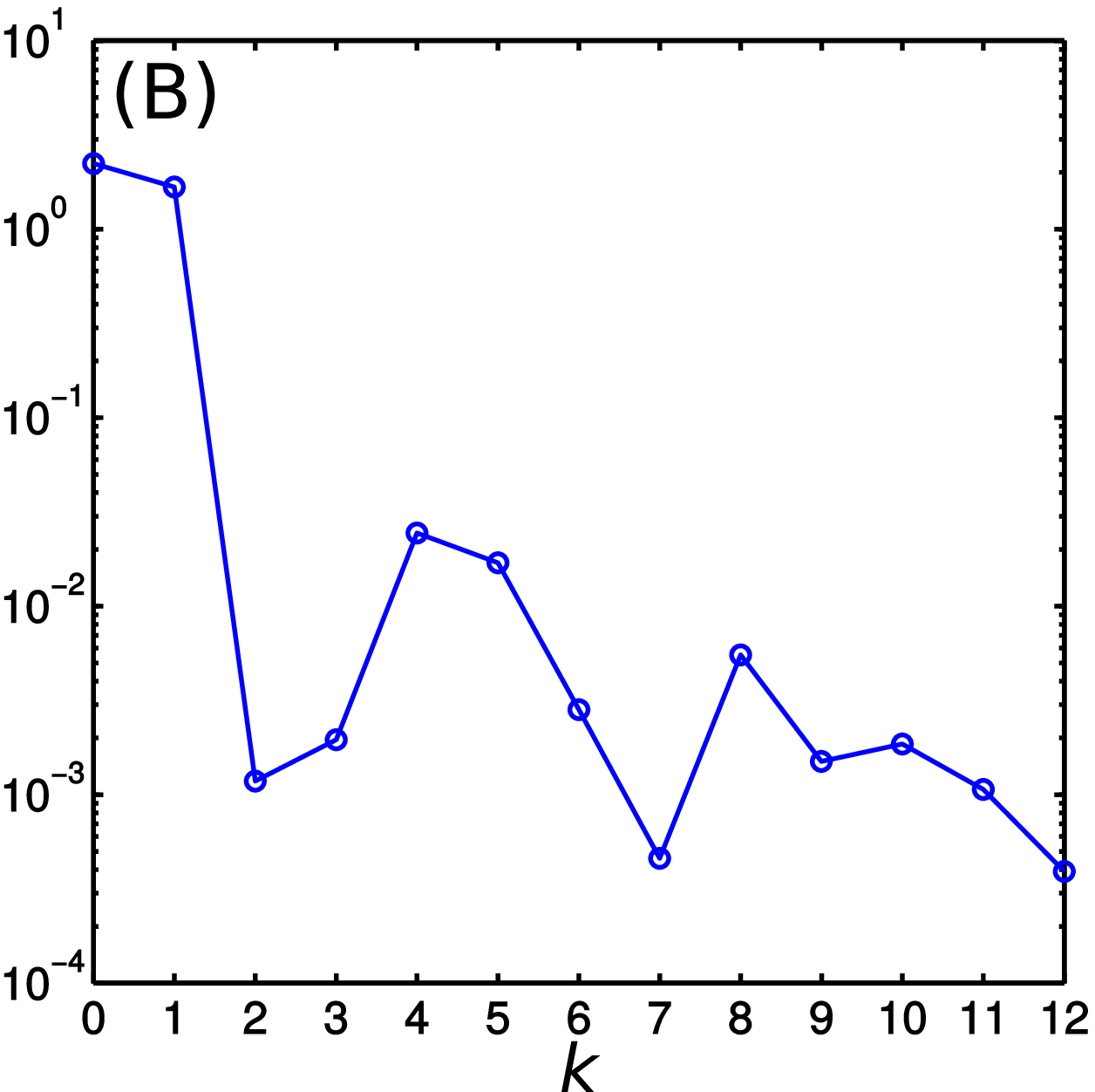} 
\par\end{centering}

\caption{A Lotka-Voltera oscillator with slowly varying frequency and amplitude, Example~\ref{sec:VL_oscillator}.
  (A) The level sets of the slow variable $I$ projected onto $x$-$y$ plane.
(B)  The absolute error in the state variables as a function of iteration.
}  \label{fig:LV}
\end{figure}

Consider the following example of a Voltera-Lotka oscillator
with slowly varying frequency and amplitude
\begin{equation}
\begin{aligned}
   & \dot{x} = \ei x (1-z y)  \\
   & \dot{y} = \ei z y (x-1) \\
   & \dot{z} = 0.2 x
\nonumber  
\end{aligned}
\end{equation}
Initial conditions are $(x,y,z)(0)=(1,2.9,1)$.
For fixed $z$, $(x,y)$ is a Voltera-Lotka oscillator 
whose period is of order $\ep$. 
The period and amplitude of $(x,y)$ depend on a parameter $z$,
which is given by the time integral of $x$.
As a result, $z$ is a slow variable.
It is easily verified that the first integral of the oscillator is also slow,
\begin{equation}
   I = x - \log(x) + y - \log(y)/z 
\nonumber
\end{equation}
Again, we stress that the slow variables are only used in order to demonstrate the
results of the method. 
They are {\em not} used in the numerical approximation.
In addition, Figure \ref{fig:LV}A shows the level set of the slow variable,
$\left\{ u\in\mathbb{R}^{3}:I(u)=I(x(t_{n}),y(t_{n}),z(t_{n}))\right\} $, projected onto x-y plane.
In contrast to the previous examples, the level set of the slow variable $I$ is not a circle.
As a result, $J(t)$ may have several local minima and we need to find the
first local minima which is close to the global minimum of $J$
within a few periods.
Parameters are given in Table \ref{tab:parareal-parameter-VL}.
The fine integrator is ODE45 method, and 
the coarse integrator is the Poincar\'e 2nd order multiscale method.
A $C^3$ kernel with $p=1$ is used for the filtered equation.

\begin{table}[H]
\centering{}\caption{Parareal parameters in Example~\ref{sec:VL_oscillator}.}  \label{tab:parareal-parameter-VL}
\begin{tabular}{|c|c|c|c|c|c|c|}
\hline 
$\epsilon$ & $T$ & $H$ & $h_{fine}$ & $\eta_{Poincare}$ & $h_{Poincare}$ & RelTol, AbsTol(ODE45 parameters)\tabularnewline
\hline 
\hline 
$10^{-3}$ & 10 & 1/2 & $\ep/200$ & 30$\epsilon$ & $\epsilon/10$ & $10^{-13}$, $10^{-10}$\tabularnewline
\hline 
\end{tabular} 
\end{table}

%

\subsection{Passage through resonance}
\label{sec:resonance}
One of the fundamental assumptions underling 
multiscale approaches such as Poincar\'e and other methods
is a spectral gap in the spectrum of the Jacobian of the equations of motion.
When this assumption fails, for example due to a temporary passage through resonance,
the assumption~\ref{eq:doesNotHold} may not hold close the resonance 
and typical multiscale methods fail.
However,  the applicability of the multiscale algorithms can be extended 
by the parareal approach described above, i.e., 
by resolving all scales of the dynamics - both the slow and the fast.

We consider the following example.
\begin{equation}
\begin{aligned}
   & \dot{x} = -2 \pi \ei f(z) y + 0.5 \sin (z) x\\
   & \dot{y} = 2 \pi \ei f(z) x \\
   & \dot{z} = 1
\nonumber
\end{aligned}
\end{equation}
where
 \[
    f(z) = \tanh \left( 50(z-4.5) \right).
 \]
Initial conditions are $(1,0,0)$.
In words, $f(z)$ changes smoothly from -1 to 1, vanishing at $z=4.5$.
Hence, the frequency of oscillation undergoes fast oscillations with varying 
frequency, except close to $t=4.5$.
At this time, $f(z)$ vanishes and the system is no longer highly oscillatory.
More precisely, trajectories go through a transition layer. Its width in this example 
is of order $\ep$.
The two slow variables are $I=x^2+y^2$ and $z$.
 
\begin{figure}[H]
\begin{centering}
\includegraphics[width=4.19cm]{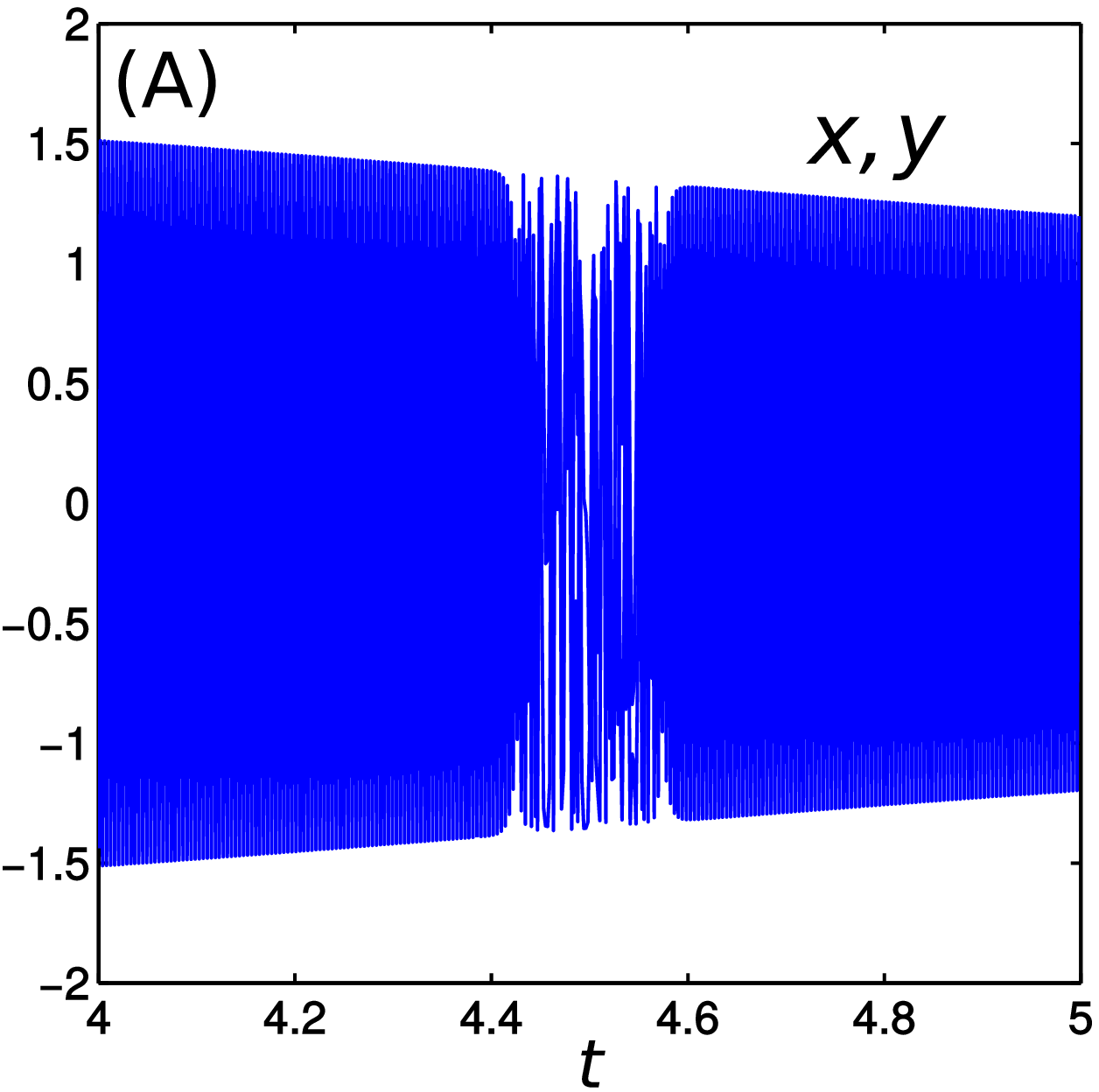} 
\includegraphics[width=4.2cm]{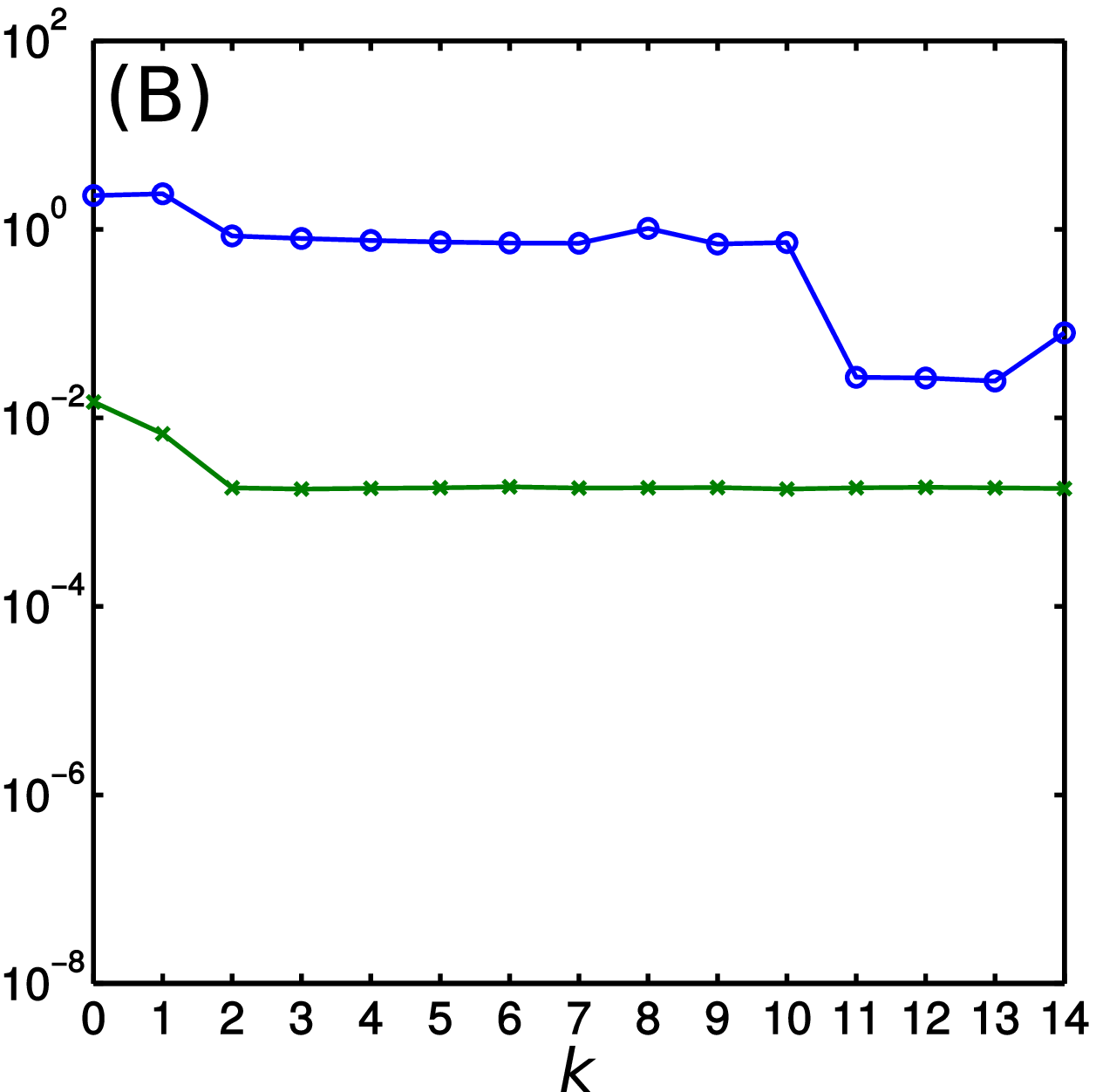}
\includegraphics[width=4.2cm]{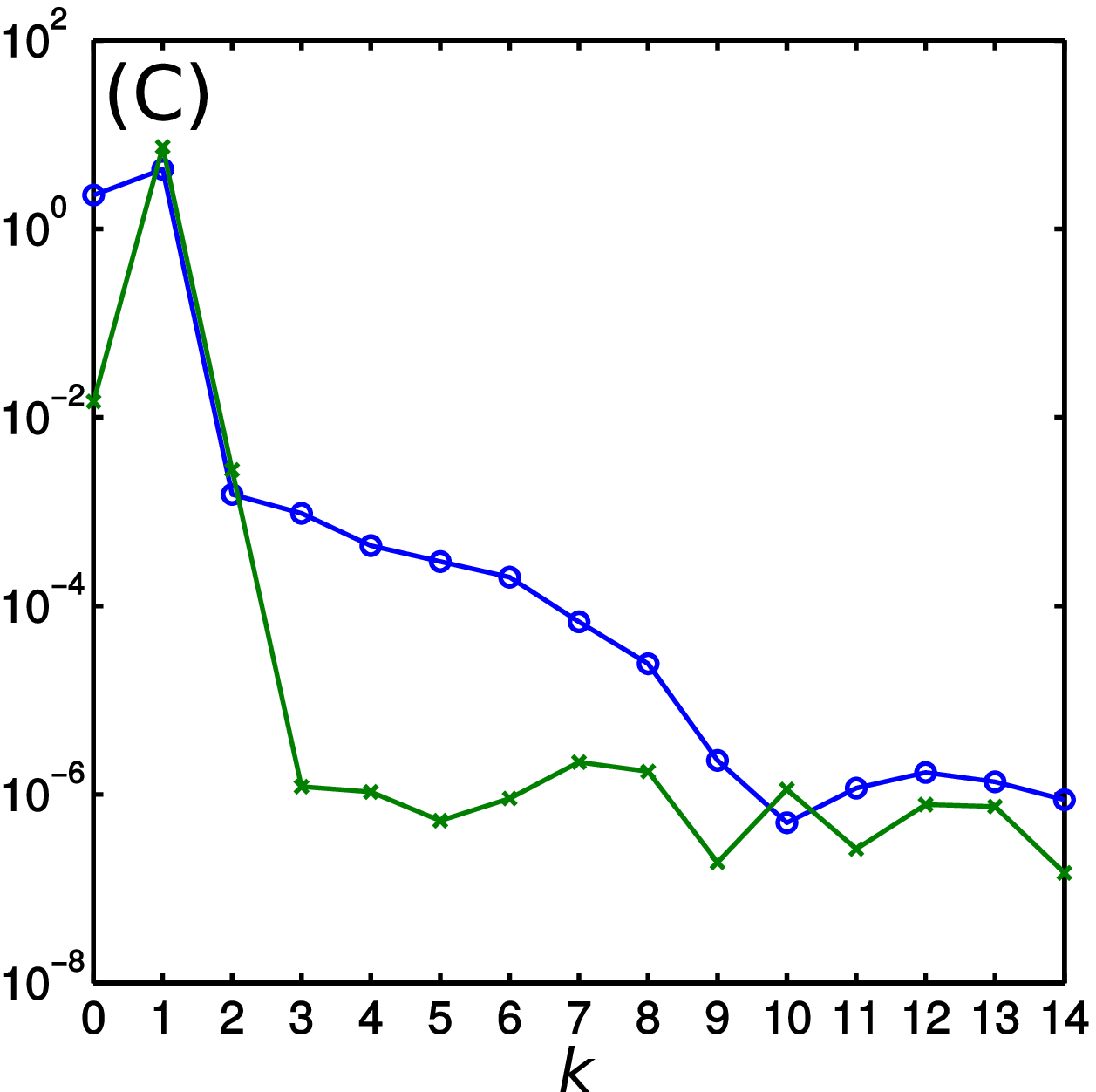}
\par\end{centering}

\caption{Passage through resonance, Example~\ref{sec:resonance}. 
(A) The solution of $x(t)$ and $y(t)$ with $t\in[4,5]$ and $\ep=10^{-4}$. The frequency function vanishes at $t=4.5$ and solutions lose their  highly oscillatory nature.
(B) The absolute errors of both state (circles) and slow (crosses) variables as a function of iterations with phase alignment at all the time.
(C) The absolute errors with phase alignment turned off near $t=4.5$. Convergence in the state variable is achieved.}
 \label{fig:layer}
\end{figure}

Figure~\ref{fig:layer}A shows the values of the state variables with $\ep=10^{-4}$.
Due to the resonance, 
the Poincar\'e method fails to capture the correct evolution of the slow variables when crossing the singular point $t=4.5$.
However, combining with parareal, the fine solution of parareal 
integrates the equation across the resonance and allows 
the multiscale method to proceed beyond the singularity. 
In Figure~\ref{fig:layer}B, the absolute error in the state variable 
does not decrease with iterations because 
the accuracy of phase alignment relies on the scale separation which does not exists near $t=4.5$.
We show, however, that the convergence in the state variable can be achieved 
with a slight modification of Algorithm~\ref{alg:pararealFull}.
Figure~\ref{fig:layer}C is obtained by skipping phase alignments, 
Step 2(c)i and ii, and replacing Step (c)iii with the naive correction \eqref{eq:pararealAlg}
in the interval near $t=4.5$.

\begin{rem}
We note that this example goes beyond the scope for which convergence
is proven in Section~\ref{seq:HiOscParareal}. 
Indeed, the purpose of the example is to demonstrate that the applicability of
the multiscale-parareal coupling using alignments may be wider than proven here.
\end{rem}

\begin{table}[H]
\centering{}\caption{Parareal parameters in  Example~\ref{sec:resonance}. \label{tab:parareal-parameter-passage}}
\begin{tabular}{|c|c|c|c|c|c|c|}
\hline 
$\epsilon$ & $T$ & $H$ & $h_{fine}$ & $\eta_{Poincare}$ & $h_{Poincare}$ & RelTol, AbsTol(ODE45 parameters)\tabularnewline
\hline 
\hline 
$10^{-4}$ & 7 & 1/4 & $\ep/200$ & 15$\epsilon$ & $\epsilon/10$ & $10^{-13}$, $10^{-11}$\tabularnewline
\hline 
\end{tabular} 
\end{table}

\subsection{The Fermi-Pasta-Ulam (FPU) problem} \label{sub:FPU}

We consider a chain of $2k$ springs on a line, connected
with alternating soft nonlinear and stiff linear springs with both ends fixed. This problem has been used as a benchmark for testing the
long-time performance of geometric integrators \cite{Hairer-Lubich-Wanner:2002}. See also \cite{Bambusi-ponno2006,Bambusi-ponno2008,Hairer-Lubich2010} and references therein for related recent work. 
The model is derived from the following Hamiltonian:
\begin{equation}
H(p,q)=\frac{1}{2}\overset{2k}{\underset{i=1}{\sum}}p_{i}^{2}+\frac{1}{4}\epsilon^{-2}\overset{k}{\underset{i=1}{\sum}}(q_{2i}-q_{2i-1})^{2}+\overset{k}{\underset{i=1}{\sum}}(q_{2i+1}-q_{2i})^{4}.
\nonumber
\end{equation}
Using the change of variables given in \cite{Ariel-Engquist-Tsai07},
the equations of motion for the system can be written as
\begin{equation}
\begin{cases}
\dot{y}_{i}=u_{i},\\
\dot{x}_{i}=\epsilon^{-1}v_{i},\\
\dot{u}_{i}=-(y_{i}-\epsilon x_{i}-y_{i-1}-\epsilon x_{i-1})^{3}+(y_{i+1}-\epsilon x_{i+1}-y_{i}-\epsilon x_{i})^{3},\\
\dot{v}_{i}=-\epsilon^{-1}x_{i}+(y_{i}-\epsilon x_{i}-y_{i-1}-\epsilon x_{i-1})^{3}+(y_{i+1}-\epsilon x_{i+1}-y_{i}-\epsilon x_{i})^{3}.
\end{cases}\label{eq:FPU}
\end{equation}
Since the ends are fixed, $y_{0}=x_{0}=y_{k+1}=x_{k+1}=0$. 

In this section, we endeavor solving \eqref{eq:FPU} in the $\mathcal{O}(\epsilon^{-1})$ time scale using the proposed methods.
The aim is to expose the limit of the various different algorithmic components of the proposed methods. 
Clearly, the time scale in which we compute the solutions of \eqref{eq:FPU} is out of the scope of the analysis that we presented earlier. 

System \eqref{eq:FPU} is solved in $\R^8$ with $k=2$, and it admits seven slow variables --- 
they are the total energies of the stiff springs,
$I_{i}=x_{i}^{2}+v_{i}^{2}$ for $i=1,2$, the relative phases between the stiff springs, $\phi=x_1 x_2 + v_1 v_2$, 
and all the degrees of freedom which are related to the soft springs: $y_i$ and $u_i$, $i=1,2$.
The nontrivial energy transfer and the relative phase take place in the very long $\epsilon^{-1}$ time scale.

\begin{figure}[H]
\begin{centering}
\includegraphics[width=5cm]{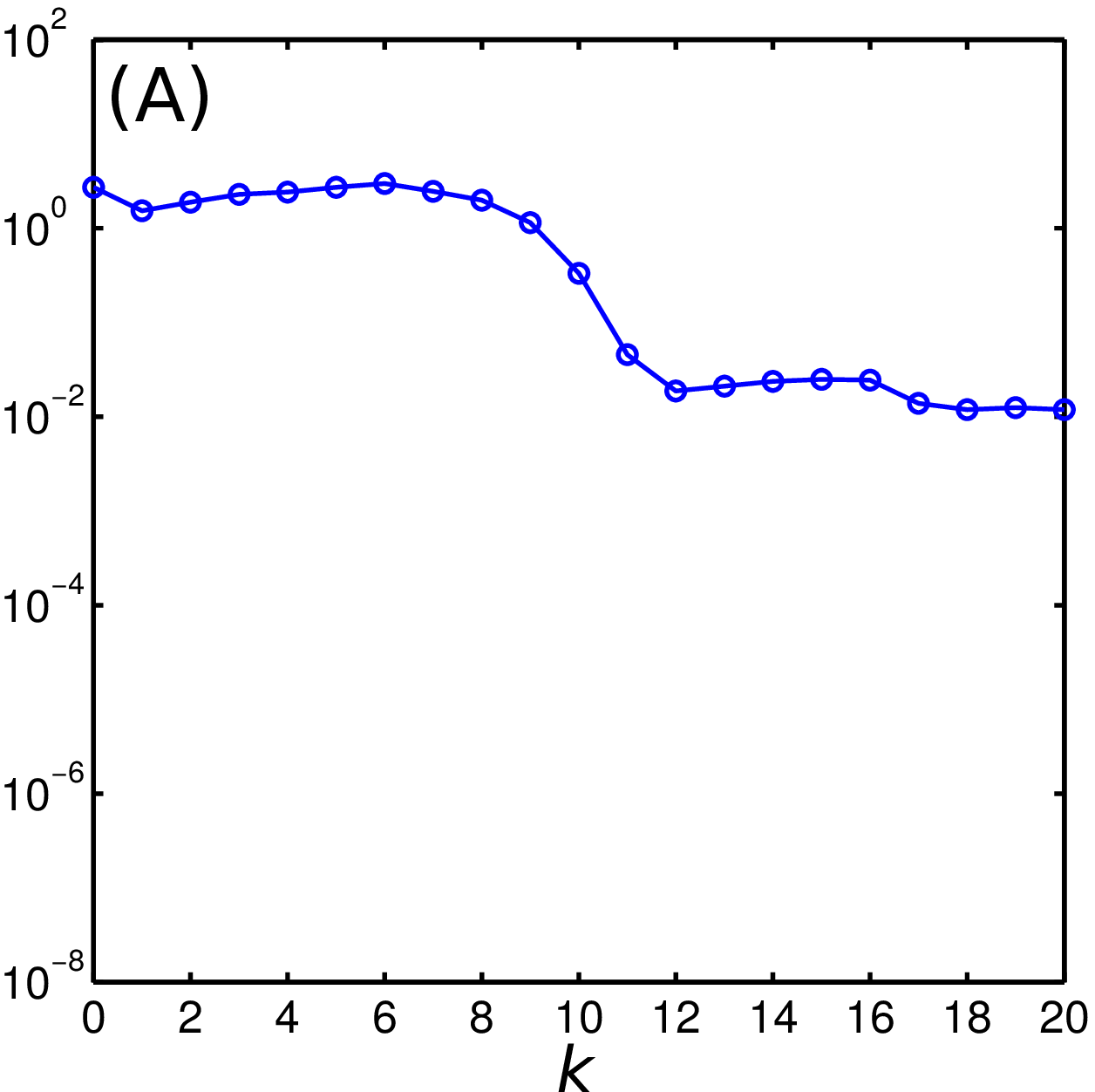}  \;
\includegraphics[width=5cm]{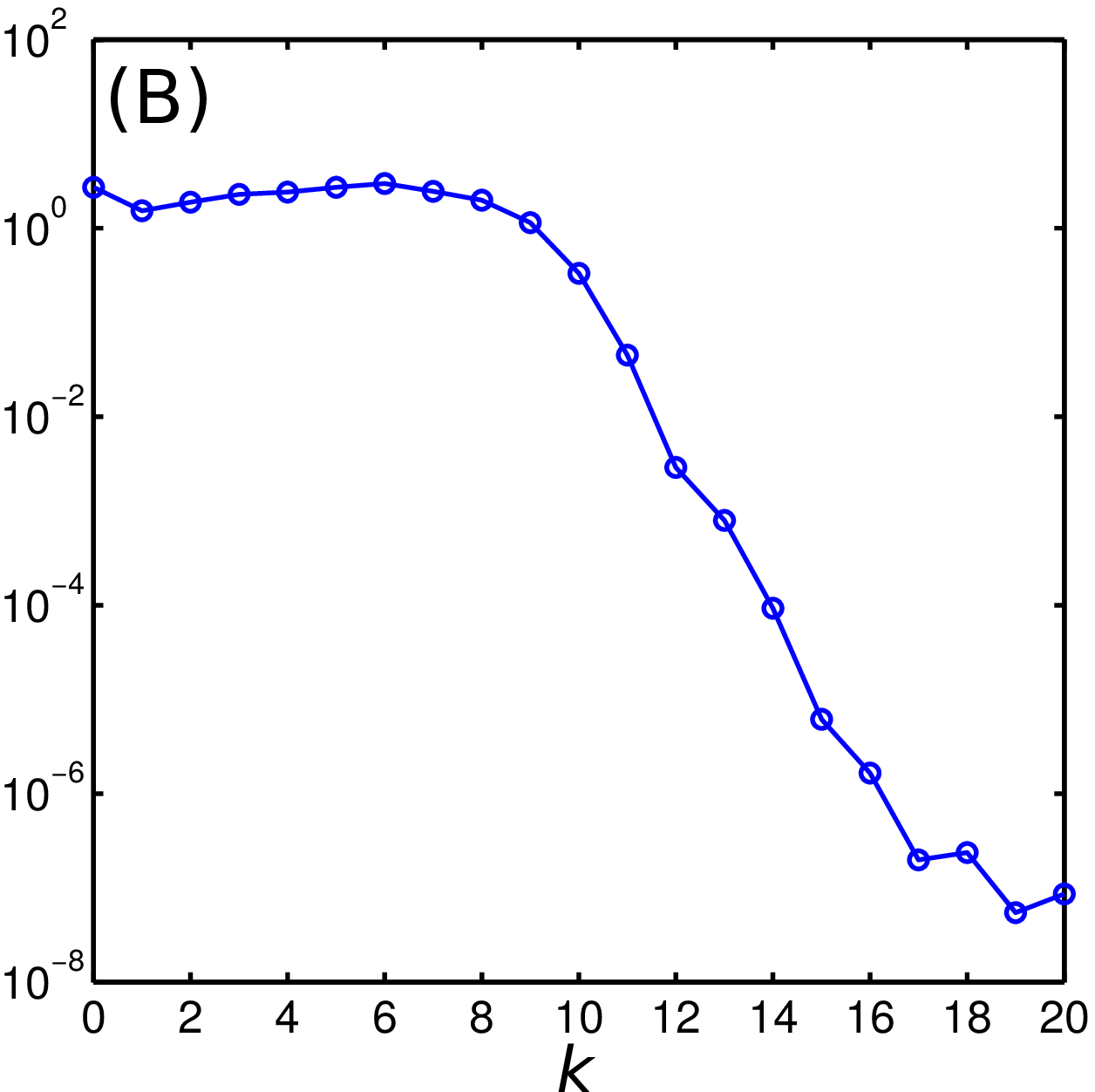} 
\par\end{centering}

\caption{\label{fig:FPU-parareal}The Fermi-Pasta-Ulam problem, Example~\ref{sub:FPU}.
The absolute error in the state variables as a function of iteration is depicted. (A) Phase alignments with the full system \eqref{eq:FPU}. 
(B) Phase alignments with the unperturbed equation explained in the paragraph.}
\end{figure}

In Figure \ref{fig:FPU-parareal}A, we present the maximum errors in the state variable 
 in a long time interval $[0,T]=[0,\epsilon^{-1}/2]$, with $\epsilon=10^{-3}$.
The system~\eqref{eq:FPU} is integrated using Algorithm~\ref{alg:pararealFull}, 
applying the global alignment algorithm described in Section~\ref{seq:HiOscParareal}
to achieve convergence in the state variable.
The results are computed with the initial conditions $y_1=x_1=y_2=x_2=1$, $(u_1,v_1,u_2,v_2)=(0,1.2,1,0)$ and with the parameters given in Table~\ref{tab:FPU}.  The Verlet method was used as the fine integrator and 
the Poincar\'e 2nd order method (Verlet macro-solver and ODE45 micro-solver) as the multiscale integrator. We further point out that the errors will decrease as $\epsilon$ becomes smaller.

In Figure \ref{fig:FPU-parareal}B, we present the maximum errors computed using the same set of conditions and parameters, \emph{except that the phase alignments are computed by the numerical solutions of the modified equation}: 
\begin{equation}
\begin{cases}
\dot{y}_{i}=0\cdot u_{i},\\
\dot{x}_{i}=\epsilon^{-1}v_{i},\\
\dot{u}_{i}=0\cdot \left\{-(y_{i}-\epsilon x_{i}-y_{i-1}-\epsilon x_{i-1})^{3}+(y_{i+1}-\epsilon x_{i+1}-y_{i}-\epsilon x_{i})^{3} \right\},\\
\dot{v}_{i}=-\epsilon^{-1}x_{i}+0 \cdot \left\{(y_{i}-\epsilon x_{i}-y_{i-1}-\epsilon x_{i-1})^{3}+(y_{i+1}-\epsilon x_{i+1}-y_{i}-\epsilon x_{i})^{3}\right\}, 
\end{cases}\label{eq:FPU-0}
\end{equation}
with $y_{0}=x_{0}=y_{k+1}=x_{k+1}=0$. 
This is a slightly more general type of \lq\lq{}unperturbed\rq\rq{} equations as those defined in Section~5. This approach may be thought of as phase alignments with some slow variables constrained. A detailed discussion about such topic is out of the scope of this paper, and shall be addressed in a future work.

\begin{rem}
We observe that in this study, phase alignments by solutions of \eqref{eq:FPU-0} yield significantly superior results. Indeed,
our experience indicated that for the problems that can be solved by the Poincar\'e method described in Section~\ref{sec:reviewPoincare}, it is generally better to use the so-called \lq\lq{}unperturbed\rq\rq{} equations for phase alignments, as the slow variables are constrained and would not 
deteriorate the performance of the phase alignment steps. However, we chose to present a more general setup in this paper, avoiding any specific choice of coarse/multiscale integrators in describing the algorithms.\end{rem}


\begin{table}[H]
\centering{}\caption{\label{tab:FPU}Parareal parameters in Example \ref{sub:FPU}.}
\begin{tabular}{|c|c|c|c|c|c|c|c|}
\hline 
$\epsilon$ & $T$ & $H$ & $h_{fine}$ & $\eta_{Poincare}$ & $h_{Poincare}$ & RelTol, AbsTol(ODE45 parameters)\tabularnewline
\hline 
\hline 
$10^{-3}$ & $\epsilon^{-1}/2$ & $1/4$ & $\epsilon/20$ & $10\epsilon$ & $\epsilon/20$ & $10^{-12}$, $10^{-10}$\tabularnewline
\hline 
\end{tabular}
\end{table}

\section{Summary}
\label{sec:summary}
\setcounter{equation}{0}

The paper describes two approaches to incorporate multiscale integrators as coarse integrators in
 parareal methods.
The first, presented in Section~\ref{sec:onlySlow}, approximates all the slow
variables.
However, the numerical approximation of the state variables $u_n^k$ does not 
converge to the true solution $u(nH)$.
This parareal-multiscale combination has several advantages
compared to other multiscale schemes.
\begin{itemize}
\item It offers increased stability and is less sensitive to
  the choice of parameters. 
  Intuitively, the parareal iterations can "fix" errors 
  incurred by the inexact multiscale scheme.
\item It offers increased accuracy.
   In fact, the accuracy of slow variables may be smaller than ${\mathcal O}(\ep)$,
     which is a theoretical limit for Poincar\'e and other multiscale methods
     that are based on averaging or homogenization principles.
\item It may be applied for systems with moderate scale separation. 
   Most multiscale methods are more efficient than conventional, non-multiscale
   schemes if the separation in scale is large enough, i.e., if $\ep$ is sufficiently small.
   However, they typically become less efficient 
   or unstable at intermediate values of $\ep$.
\item The method may be used in situations in which the dynamics looses
   its multiscale structure in a short transition layer, for example, 
   due to passage through resonance,
   see Example \ref{sec:resonance}.
\end{itemize}

The second approach, presented in Section~\ref{eq:phaseCont}, computes convergent approximation 
to all state variables in the system.
This algorithm requires the phase alignment procedure, described in Section~\ref{seq:HiOscParareal} in addition to the steps needed in our first algorithm.
We prove that the accuracy of the scheme in the sup norm after $K$ iterations is of order
$\ep^{-1} H^K+E_f$, where $H$ is the a coarse step size.
In particular, the number of iterations to achieve a given error tolerance 
is logarithmic in $\ei$.

The computational cost of the method 
can be divided into two contributions.
The first is the cost of the fine integrator invoked at each parareal iteration.
With $K$ iterations its contribution to the overall
cost is proportional to $K H \ei$.
The second contribution comes from the overhead of coarse multiscale integrators
and phase alignment.
While this contribution is independent of $\ep$, it grows linearly with
the number of coarse step sizes, $H^{-1}$.
Hence, there is a trade off in choosing $H$.
With a large scale separation $\ep \ll 1$, the first contribution dominates
and, assuming maximal parallelization is available, it is advantageous to use a relatively small $H$,
even if the multiscale method allows larger steps.
The two contributions balance if one takes $H=\sqrt{\epsilon}$, which implies a computational
cost of order $K\epsilon^{-1/2}$. 
In contrast, parallel methods using conventional integrators will require at least $\mathcal{O}(\epsilon^{-1})$ steps.

The multiscale-parareal coupling strategy using alignments is not limited to the proposed particular implementation involving a single-frequency fast variable. 
In the multi-frequency case, different alignment algorithms may be possible. For example, if the slow-fast system is given explicitly, one could simply set the fast variable in the coarse integrator to coincide with the fine one.


\section*{Acknowlegments}
Tsai's research is  supported by a Simons Foundation Fellowship, NSF grants DMS-1217203, and DMS-1318975.

%

\appendix

\section{Appendix: Adaptive search algorithm} \label{appen:search_algorithm}
\setcounter{equation}{0}

We explain our implementation for the local alignment ${\mathcal{S}}_{0}$. 
This algorithm adaptively searches local minima of the functional $J(t)$ by 
adjusting the step size and the computational interval.
The first part of the algorithm numerically solves the $l_2$ minimization problem of $J(t)$
using quadratic interpolation. The second part
computes a convex combination for the local alignment ${\mathcal{S}}_{0}(u_{0};v_{0})$.

For simplicity,
we assume that the fine solver applies a numerical method with step
size $h={\mathcal{O}}(\epsilon)$. Also, let $\eta=\eta_{phase}>0$ denote a parameter
that is larger than at least one period. It is assumed that $\eta$
is of order $\epsilon$, but the size of $\eta$ is not explicitly known. 
Our algorithm adaptively finds the size of $\eta$ and
locally integrates the HiOsc ODE in short time segments of length
$[-\eta,\eta]$. Thus, the computational cost of each alignment procedure
is independent of~$\epsilon$.

Recall that at $t=0$, given two points $u_{0}$ and $v_0$ such that
$|\xi(u_0) - \xi(v_0)| = {\mathcal O} (\ep)$,
local alignment approximates a new point  
$w_{0}$ such that $\xi(w_{0})=\xi(u_{0})$ and $\phi(w_{0})=\phi(v_{0})$.
In addition, first forward and backward alignment times $t_{0}^{\pm}$  are sought,
in which $| {\mathcal F}_{t_0^\pm} u_0 - v_0|$ are local minima. 

In the following, denote by $[x]$ the integer value of $x$. 
\begin{algorithm}
Local alignment ${\mathcal{S}}_{0}(u_{0};v_{0})$
\end{algorithm}
\begin{enumerate}
\item Adaptive search for first two local minima of $J(t)$ associated with positive(+) and negative(-) orientations.
 \begin{enumerate}
	\item Set $\eta=\epsilon$, $h=\epsilon/100$, $Tol=\infty$, $u_{prev}=\infty$,
	and $I_{\pm}=\{\}.$
	\item While $Tol>\ep/100$ or $I_{\pm}\neq\{\}$
	\begin{enumerate}
	\item Forward/Backward local integration: Compute $u_{0,i}=\mathcal{F}_{ih}u_{0}$,
	$i\in \{ 0, \dots , [\pm\eta/h] \}.$ 
	\item Compute the local minimum of $J(i)=|u_{0,i}-v_{0}|^{2},$ $i \in \{ 0, \dots , [\pm\eta/h] \}$.
	\\
	Denote it $I_{\pm}=(i_{1\pm},i_{2\pm},\cdots,i_{r\pm})$. 
	\item If $I_{\pm}\neq\{\}$, set $Tol=|u_{prev}-u_{0,i_{1\pm}}|$, $u_{prev}=u_{0,i_{1\pm}},$ and $h=h/2$.\\
	Else, set $\eta=2\eta$ and $h=\ep/100$.
	\end{enumerate}
	End while.
	\item Quadratic interpolation: For each of the first two indices $i_{1\pm}$, improve
	the minimization of $J$ using a quadratic interpolation. Let $p_{\pm}(t)$
	denote the polynomial such that 
	\[
	p_{\pm}(i_{1\pm}+j)=u_{0,i_{1\pm}+j},\quad j=-1,0,1.
	\]
Denote by $t_0^{\pm}$ the minimum of $|p_{\pm}(t)-v_{0}|$.
\end{enumerate}
%
%
	\item Let 
	\[
	\lambda_{+}=\frac{t_{0}^{-}}{t_{0}^{+}+t_{0}^{-}}~~,~\lambda_{-}=\frac{t_{0}^{+}}{t_{0}^{+}+t_{0}^{-}},
	\]
	and define
	\begin{equation}
	{\mathcal S}_0 (u_0;v_0) = \lambda_+ {\mathcal F}_{t_0^+} u_0 + \lambda_- {\mathcal F}_{t_0^-} u_0.
	\label{eq:adaptive_s0}
	\end{equation}
 \end{enumerate}
\noindent
In \eqref{eq:adaptive_s0}, we implement the fine integrator ${\mathcal F}$ with step size $h$ determined in step 1(b). 
If $t_0^\pm/h$ is not integer, then ${\mathcal F}$ is to be evaluated using quadratic interpolation.

\section{Appendix: Convergence of the symmetric Poincar\'e method}
\setcounter{equation}{0}

To prove convergence of the symmetric Poincar\'e method algorithm described in Section~\ref{sec:z},  we use a diffeomorphism  $\Psi : u \rightarrow (\xi(u),\phi(u))$ 
given in \eqref{eq:odefastslow},
\begin{equation}
   \begin{cases}
      \dot{\xi}=g_{0}(\xi,\phi), & \xi(0)=\xi(u_0),\\
      \dot{\phi}=\epsilon^{-1}g_{1}(\xi)+g_{2}(\xi,\phi), & \phi(0)=\phi(u_0),
   \end{cases} 
\label{eq:slow-fast-app}
\end{equation}
where $g_0 (\xi,\phi)$ and $g_2 (\xi,\phi)$ are 1-periodic in $\phi$.
We stress that the variables $(\xi,\phi)$ are only used in the analysis but {\em not} in the numerical algorithm.
In Section~\ref{sec:hiosc}, the dynamics of the slow variables can be approximated by an averaged equation of 
\eqref{eq:averagedEqn},
\begin{equation}
\begin{aligned}
   & \dot{\bar \xi} = F(\bar \xi),~~~ F({\bar \xi}) = \int g_0(\xi,\phi) d \phi_\xi,\\
   & \bar\xi(0)=\xi(u_0).
\end{aligned}
\label{eq:avr-eqn-app}
\end{equation}
Recall that the construction given in Section~\ref{sec:z},
\begin{equation}
   u_{n+1} = \gamma_{-1}^* + \frac{H}{2\eta} \left( \gamma_1^* - \gamma_{-1}^* \right).
\nonumber
\end{equation}
%
To simplify the calculation, we generate $\gamma^*$ using the symmetric shape ($z$-shape) which is centered at $\gamma^*_{0}$. See Figure~\ref{fig:BFHMM_schematics}.
The formulas for $\gamma^*_{-1}$, $\gamma^*_{0}$ and $\gamma^*_{1}$ are thus of the form
\begin{equation}
   \gamma_{-1}^*= {\mathcal F}^0_{\eta} {\mathcal F}_{ -\eta} u_n,~~
   \gamma_{0}^* = u_n,~~
   \gamma_{1}^* = {\mathcal F}^0_{-\eta} {\mathcal F}_{ \eta} u_n .
\nonumber
\end{equation}
We then prove the following theorem.
\begin{thm}
The Poincar\'e force estimator defined by 
\[
\mathcal{P} =  \frac{1}{2\eta} \left( {\mathcal F}^0_{-\eta} {\mathcal F}_{\eta}
-  {\mathcal F}^0_{\eta}  {\mathcal F}_{-\eta} \right)
\]
satisfies the following estimates in the $(\xi,\phi)$ coordinate:
\begin{equation}
\xi(\mathcal{P}u_n)  = \frac{d}{dt} \bar\xi (t_n) +E_F,
~~~\left| \phi(\mathcal{P}u_n)  \right|=\mathcal{O}\left(\frac{\eta^{2}}{\epsilon}\right)
\label{eq:thm-z}
\end{equation}
where $E_F$ denotes the error of the filtered equation in approximating an averaged equation.
\label{thm:z-shape}
\end{thm}

\noindent

\begin{proof}
Consider the unperturbed system associated with  \eqref{eq:slow-fast-app},
\begin{equation}
\begin{cases}
\dot{\xi}=0, & \xi(0)=\xi_{0},\\
\dot{\phi}=\frac{1}{\epsilon}g_{1}(\xi), & \phi(0)=\phi_{0}.
\end{cases}\label{eq:appen_unper}
\end{equation}
%
When $\mathcal{F}$ is corresponding to the filtered equation \eqref{eq:filtered-eqn}, it is proved in \cite{BFHMM2012} that
$\mathcal{F}_{\pm\eta} u_n$ is essentially close to averaged $\bar{\xi}(t_n \pm \eta)$, respectively.
We denote by $(\xi^+,\phi^+)$  the slow-fast coordinate of the point $\gamma_{1}^{*}$, i.e., $(\xi^+,\phi^+) = (\xi(\gamma^*_1),\phi(\gamma^*_1))$. Similarly, $(\xi^-,\phi^-)=(\xi(\gamma^*_{-1}),\phi(\gamma^*_{-1}))$.
Then $\xi^{+}$ and $\phi^{+}$ satisfy
\[
\begin{cases}
\xi^{+}=\xi_{0}+\int_{0}^{\eta}F(\bar{\xi}^{F}(t))dt,\\
\phi^{+}=\phi_{0}+\int_{0}^{\eta}\frac{1}{\epsilon}g_{1}(\bar{\xi}^{F}(t))dt+\int_{0}^{\eta}G(\bar{\xi}^{F}(t))dt-\frac{\eta}{\epsilon}g_{1}(\xi^{+}),
\end{cases}
\]
where $\bar{\xi}^{F}(t)$ corresponds to the solution of \eqref{eq:avr-eqn-app} forward in time and $G$ is the averaged $g_2(\xi,\phi)$ due to the filtered equation.
On the other hand, for $\gamma_{-1}^{*}$,
\[
\begin{cases}
\xi^{-}=\xi_{0}-\int_{0}^{\eta}F(\bar{\xi}^{B}(t))dt,\\
\phi^{-}=\phi_{0}-\int_{0}^{\eta}\frac{1}{\epsilon}g_{1}(\bar{\xi}^{B}(t))dt-\int_{0}^{\eta}G(\bar{\xi}^{B}(t))dt+\frac{\eta}{\epsilon}g_{1}(\xi^{-}),
\end{cases}
\]
where $\bar{\xi}^{B}(t)$ corresponds to the backward in time solution of \eqref{eq:avr-eqn-app}.

For the slow variables, evaluating the force by $\gamma^*_1 - \gamma^*_{-1}$ yields
\[
\xi^{+}-\xi^{-}=\int_{0}^{\eta}F(\bar{\xi}^{F}(t))dt+\int_{0}^{\eta}F(\bar{\xi}^{B}(t))dt
\]
which approximates the evolution of $\bar{\xi}$ over the interval $[-\eta,\eta]$.
For the fast variable,
\begin{eqnarray*}
\phi^{+}-\phi^{-} & = & \underset{I_1}{\underbrace{\int_{0}^{\eta}\frac{1}{\epsilon}g_{1}(\bar{\xi}^{F}(t))dt-\frac{\eta}{\epsilon}g_{1}(\xi^{+})+\int_{0}^{\eta}\frac{1}{\epsilon}g_{1}(\bar{\xi}^{B}(t))dt-\frac{\eta}{\epsilon}g_{1}(\xi^{-})}}\\
&  &+ \underset{I_2}{\underbrace{\int_{0}^{\eta}G(\bar{\xi}^{F}(t))dt-
\int_{0}^{\eta}G(\bar{\xi}^{B}(t))dt}}.
\end{eqnarray*}
By considering $\bar{\xi}^{F}(t)=\xi_{0}+\int_{0}^{t}F(\bar{\xi}^{F})ds$ and
$\bar{\xi}^{B}(t)=\xi_{0}-\int_{0}^{t}F(\bar{\xi}^{B})ds$,
we can compute the Taylor series to estimate $\epsilon I_1$,
developed at $\xi_0$.

\begin{eqnarray*}
\epsilon I_1&=&\left\{ \int_{0}^{\eta}g_{1}\left(\xi_{0}+\int_{0}^{t}F(\bar{\xi}^{F}(s))ds\right)dt-\eta g_{1}\left(\xi_{0}+\int_{0}^{\eta}F(\bar{\xi}^{F}(t))dt\right)\right\} \\
&&+\left\{ \int_{0}^{\eta}g_{1}\left(\xi_{0}-\int_{0}^{t}F(\bar{\xi}^{B}(s))ds\right)dt-\eta g_{1}\left(\xi_{0}-\int_{0}^{\eta}F(\bar{\xi}^{B}(t))dt\right)\right\} \\
&=& \left\{\int_{0}^{\eta}\left(\int_{0}^{t}\nabla{g}_{1}(\xi_{0})^T F(\bar{\xi}^{F}(s))ds\right)dt-\eta\int_{0}^{\eta}\nabla{g}_{1}(\xi_{0})^T F(\bar{\xi}^{F}(t))dt\right\} \\
&& -\left\{ \int_{0}^{\eta}\left(\int_{0}^{t}\nabla{g}_{1}(\xi_{0})^T F(\bar{\xi}^{B}(s))ds\right)dt-\eta\int_{0}^{\eta}\nabla{g}_{1}(\xi_{0})^T F(\bar{\xi}^{B}(t))dt \right\} +\mathcal{O}(\eta^{3}).
\end{eqnarray*}
Here, the terms $\eta g_{1}\left({\xi}_{0}\right)$ canceled. Due to the symmetric structure, expanding both $\bar{\xi}^{F}(s)$ and $\bar{\xi}^{B}(s)$ about $s=0$ using Taylor series, all $\mathcal{O}(\eta^2)$ terms vanish and we have

\[
\epsilon I_1=\mathcal{O}(\eta^{3}).
\]
Similarly, one can show that
\[
I_2=\mathcal{O}(\eta^{2}).
\]
Therefore, there exist nonnegative constants $C_1$ and $C_2$ such that
\begin{equation}
\left|\phi(\gamma_{1}^{*})-\phi(\gamma_{-1}^{*})\right|=|\phi^+ - \phi^- | \leq C_1 \frac{\eta^{3}}{\epsilon}+C_2 \eta^2.
\label{eq:2ndestimate}
\end{equation}
This completes the proof. 
\end{proof}
\noindent
It is important to
use an appropriate sequence of interpolating points in the state space, $\gamma_{k}^{*},$ as it directly impacts on the accuracy of Poincar\'e method.
With a parameter $\eta=C\epsilon$, \eqref{eq:thm-z} shows that the
force estimator $\mathcal{P}$
generates two interpolating points to estimate the evolution of the slow variables over $2\eta$ with an $\mathcal{O} (\epsilon)$ disagreement in the fast variable.
This results in stable and accurate approximations. 
However, we remark that 
the force estimation using the algorithm \eqref{basic_FE}
sometimes introduces an $\mathcal{O} (1)$ difference in the fast variable and 
thus shifts the slow variables.

\bibliographystyle{plain}
\addcontentsline{toc}{section}{\refname}\bibliography{odes}

%

%
%


\end{document}